\documentclass[11pt]{article}
\usepackage{amssymb, amsthm, amsmath, amscd}

\newtheorem{thm}{Theorem}[section]
\newtheorem{lem}[thm]{Lemma}
\newtheorem{prop}[thm]{Proposition}
\newtheorem{cor}[thm]{Corollary}
\theoremstyle{definition}\newtheorem{df}[thm]{Definition}
\theoremstyle{definition}\newtheorem{rem}[thm]{Remark}
\theoremstyle{definition}\newtheorem{exm}[thm]{Example}

\renewcommand{\phi}{\varphi}

\newcommand{\N}{\mathbb{N}}
\newcommand{\Z}{\mathbb{Z}}
\newcommand{\Q}{\mathbb{Q}}
\newcommand{\R}{\mathbb{R}}
\newcommand{\T}{\mathbb{T}}

\newcommand{\Homeo}{\operatorname{Homeo}}
\newcommand{\Aut}{\operatorname{Aut}}
\newcommand{\Aff}{\operatorname{Aff}}
\newcommand{\Isom}{\operatorname{Isom}}
\newcommand{\id}{\operatorname{id}}
\newcommand{\Coker}{\operatorname{Coker}}

\newcommand{\dist}{\operatorname{dist}}
\newcommand{\cel}{\operatorname{cel}}
\newcommand{\diag}{\operatorname{diag}}
\newcommand{\sdp}{\operatorname{sdp}}
\newcommand{\SP}{\operatorname{SP}}

\newcommand{\dt}{\delta}
\newcommand{\ep}{\varepsilon}
\newcommand{\xa}{(X,\alpha)}
\newcommand{\yb}{(Y,\beta)}
\newcommand{\CA}{$C^*$-algebra}
\newcommand{\SCA}{$C^*$-subalgebra}
\newcommand{\hm}{homomorphism}
\newcommand{\morp}{completely positive linear map}

\title{Minimal dynamical systems on the product \\
of the Cantor set and the circle II}
\author{Huaxin Lin and Hiroki Matui}
\date{}

\begin{document}
\maketitle

\begin{abstract}
Let $X$ be the Cantor set and $\phi$ be a minimal homeomorphism
on $X\times\T$.
We show that the crossed product $C^*$-algebra $C^*(X\times\T,\phi)$ is
a simple $A\T$-algebra provided that the associated cocycle takes
its values in rotations on $\T$.
Given two minimal systems $(X\times\T,\phi)$ and $(Y\times\T,\psi)$
such that $\phi$ and $\psi$ arise from cocycles with values
in isometric homeomorphisms on $\T$,
we show that two systems are approximately $K$-conjugate
when they have the same $K$-theoretical information.
\end{abstract}

\section{Introduction}

It has been known that the study of minimal topological dynamical systems
is related to the study of the associated simple crossed product \CA s.
Indeed, J. Tomiyama \cite{Tm} proved
that, if $(X,\alpha)$ and $(Y,\beta)$ are
two topological transitive dynamical systems, then
they are flip conjugate if and only if there is an isomorphism between
the crossed product \CA s which maps $C(X)$ onto $C(Y).$
With the development of the classification of simple amenable \CA s,
it becomes possible to have some $K$-theoretical description
of some interesting equivalence relation among minimal dynamical systems.
In fact, Giordano, Putnam and Skau in \cite{GPS},
in the case that $X$ and $Y$ are Cantor sets, among other things,
showed that strong orbit equivalence can be determined by $K$-theory
of the dynamical systems.
In \cite{LM1}, we showed that, for Cantor minimal systems,
the strong orbit equivalence is equivalent to the approximate $K$-conjugate.
Both results used the fact that the crossed product \CA s arising
from Cantor minimal systems are simple $A\T$-algebras with real rank zero.
It seems that the notion of approximate $K$-conjugacy is not only
closely related to the above mentioned result of Tomiyama but also
closely related to that of Giordano, Putnam and Skau.
Moreover, it seems possible that, for more general spaces,
at least in connection with \CA\, theory,
versions of approximate conjugacy may be more interesting relations
than that of conjugacy or even strong orbit equivalence.
It seems also possible that, for example, approximate $K$-conjugacy
may be determined by the $K$-theoretical data of the dynamical systems
in much more general situation.
As a preliminary attempt, in \cite{LM2}, we studied
the minimal dynamical systems $(Y,h),$ where $Y=X\times \T.$
Since the Cantor set is totally disconnected and $\T$ is connected,
$h$ can be written as $h=\sigma\times\phi,$
where $\sigma$ is a minimal homeomorphism on $X$ and
$\phi_x$ is a homeomorphism on $\T$ for each $x\in X.$
We showed in \cite{LM2} that $K$-theoretical data of the minimal systems
determines the approximate $K$-conjugacy in the case that
$h$ is rigid and $\phi_x$ is a rotation for each $x\in X.$

In this paper, we first consider the case that $h=\alpha\times R_{\xi}$
($\xi\in C(X,\T)$) which may not be rigid.
We show that the crossed products have tracial rank no more than one.
Consequently, they are simple $A\T$-algebras.
One of problems related to the proof is to answer the following question:
Let $u_1$ and $u_2$ be two unitaries in a unital simple separable \CA\,
$A$ with tracial rank no more than one.
When are they approximately unitarily equivalent?
In the case that $A$ has tracial rank zero, it is known
that $u_1$ and $u_2$ are approximately unitarily equivalent
if and only if $[u_1]=[u_2]$ in $K_1(A)$ and
$\tau\circ f(u_1)=\tau\circ f(u_2)$ for all continuous functions
$f\in C(S^1)$ and all tracial states $\tau\in A$.
Let $CU(A)$ be the closure of the commutator subgroup of $U(A).$
If $u_1$ and $u_2$ are approximately unitarily equivalent,
then $\overline{u_1}=\overline{u_2}$ in $U(A)/CU(A).$
In the case that $TR(A)=0,$ $U_0(A)/CU(A)=\{0\}.$
However, when $TR(A)=1,$ this is no longer the case.
We prove that in this case,
$u_1$ and $u_2$ are approximately unitarily equivalent
if and only if $[u_1]=[u_2]$ in $K_1(A),$
$\tau\circ f(u_1)=\tau\circ f(u_2)$ for all continuous functions
$f\in C(S^1)$ and all tracial states $\tau\in A$, and
$\overline{u_1}=\overline{u_2}$ in $U(A)/CU(A).$

The rest of this paper studies the problem
when two minimal systems $(X\times\T,\alpha\times\phi)$ and
$(Y\times\T,\beta\times\psi)$ are approximately $K$-conjugate.
Roughly speaking, these two systems are approximately $K$-conjugate,
if there are two sequences of homeomorphisms
$\sigma_n:X\times\T\to Y\times\T$ and $\gamma_n:Y\times\T\to X\times\T$
such that
\[ \lim_{n\to\infty}\lVert f\circ\sigma_n\circ\alpha\circ\sigma_n^{-1}
-f\circ\beta\rVert=0
\text{ and }
\lim_{n\to\infty}\lVert g\circ\gamma_n\circ\beta\circ\gamma_n^{-1}
-g\circ\alpha\|=0 \]
for all $f\in C(Y\times\T)$ and $g\in C(X\times\T),$
and $\{\sigma_n\}_n$ and $\{\gamma_n\}_n$ give consistent information
on $K$-theory.
We will give a $K$-theoretical description of approximate $K$-conjugacy
in the case that $\phi_x$ and $\psi_y$ are isometries on $\T.$
We will apply some results and methods
in the theory of classification of simple amenable \CA s.

It was shown in \cite{LM2} that $(X\times\T,\alpha\times\phi)$ is rigid
if and only if the corresponding crossed product has real rank zero.
As an application of a result of N. C. Phillips, we present a proof
that the crossed product in this case actually has tracial rank zero,
whenever $\phi_x$ is an isometry on $\T$ for every $x\in X$.
Thus, simple crossed products coming from these minimal rigid systems are
covered by the classification theorem of \cite{L3}.
This paves the way to have a $K$-theoretical description of
approximate $K$-conjugacy for those minimal dynamical systems.
\bigskip

\noindent
{\bf Acknowledgements}
The first named author would like to acknowledge the support from
NSF and support from Shanghai Priority Disciplines during his
visit to East China Normal University in the summer 2004.
The second named author would like to acknowledge the support by
Grant-in-Aid for Young Scientists (B) of
Japan Society for the Promotion of Science.

\section{Preliminaries}

\begin{df}\label{Dcomu}
Let $A$ be a unital \CA. The closure of the commutator group in
$U(A)$ will be denoted by $CU(A).$
\end{df}

\begin{df}
Let $A$ and $B$ be unital \CA s and $h: A\to B$ be a unital \hm.
We will denote by $h^{\sharp}: U(A)/CU(A)\to U(B)/CU(B)$ the \hm\,
induced by $h.$
\end{df}

\begin{df}
Let $A$ be a stably finite \CA. Denote by $T(A)$ the tracial state
space. We will denote by $\Aff(T(A))$ the space of all (real) affine
continuous functions on $T(A).$

We will denote by $h_{\natural}:\Aff(T(A))\to\Aff(T(A))$ the affine
\hm\, induced by $h.$
\end{df}

\begin{df}\label{DD}
Let $Y$ be a compact metric space and
let $u\in B=M_l(C(Y))$ be a unitary.
We define
$$
D_B(u)=\min\{\|a\|:a\in B_{s.a.}\text{ such that }
\det(e^{ia}\cdot u)=1\}.
$$
If $B=\oplus_{i=1}^mM_{r(i)}(Y)$ and $u\in U(B),$ we define
$$
D_B(u)=\max\{D_{M_{r(i)}(C(Y))}(u):1\le i\le m\}.
$$
\end{df}

Let $B$ as above. If $u\in CU(B),$ then it is clear that
$D_B(u)=0.$

The following notation is taken from \cite{EG}.

\begin{df}\label{Dsdp}
Let $X$ and $Y$ be two compact metric spaces and let
$\phi:C(X)\to M_r(Y)$ be a \hm.
For each $y\in Y,$ define $\phi_y(f):C(X)\to M_r$
by $\phi_y(f)=\phi(f)(y)$ for $f\in C(X).$
There are rank one projections $e_1,e_2,\dots,e_l$ and
$x_1,x_2,\dots,x_l\in X$ ($l\le r$) such that
$\phi_y(f)=\sum_{i=1}^lf(x_i)e_i.$ Note that $x_i$ may be repeated.
Put $\SP\phi_y=\{x_1,x_2,\dots,x_l\}.$  Again, we count multiplicity
of each point in the spectrum.

For any $\eta>0$ and $\dt>0,$
a unital \hm\, $\phi:C(X)\to M_r(Y)$ is said to have
the property $\sdp(\eta,\dt)$ if for any $\eta$-ball
$$
O(x,\eta):=\{x'\in X:\dist(x',x)<\eta\}\subset X
$$
and any point $y\in Y,$
$$
\#(\SP\phi_y\cap O(x,\eta))>\dt\cdot\#(\SP\phi_y),
$$
counting multiplicity.

If $B=\oplus_{i=1}^mM_{r(i)}(C([0,1])),$
let $\pi_i:B\to M_{r(i)}(C([0,1])),$ $i=1,2,...,m.$
We say a \hm\, $\phi:C(X)\to B$ has the property $\sdp(\eta,\dt),$
if $\pi_i\circ\phi$ has the property $\sdp(\eta,\dt)$ for each $i.$
\end{df}
\bigskip

Throughout this paper, $X$ and $Y$ will be the Cantor set.
For $\alpha\in\Homeo(X)$ we denote the set of $\alpha$-invariant
probability measures on $X$ by $M_\alpha$.

Let $\alpha:X\to X$ be a homeomorphism and
$\phi:X\to\Homeo(\T)$ be a continuous map.
By $\alpha\times\phi$ we mean the homeomorphism on $X\times\T$
defined by $(x,t)\mapsto(\alpha(x),\phi_x(t))$.
It is easily seen that every homeomorphism on $X\times\T$
is of this form (see \cite[Lemma 2.1]{LM2}).
The continuous map $\phi$ is called a cocycle.
Moreover, if $\alpha\times\phi$ is minimal, then
$\alpha$ is also minimal, that is, $\xa$ is a Cantor minimal system.
Define a continuous map $o(\phi):X\to\Z_2$ by
\[ o(\phi)(x)=\begin{cases}0&\phi_x\text{ is orientation preserving}\\
1&\text{otherwise.}\end{cases} \]
We say $\alpha\times\phi$ is orientation preserving
when $o(\phi)$ vanishes in the $\Z_2$-values cohomology group
\[ C(X,\Z_2)/\{f-f\circ\alpha^{-1}:f\in C(X,\Z_2)\}. \]
Note that this group is canonically identified with
$K^0\xa\otimes\Z_2$.

Let $F:X\times\T\to X$ be the projection to the first coordinate.
Then $F$ is a factor map from $(X\times\T,\alpha\times\phi)$ to $\xa$.
We say that $\alpha\times\phi$ is rigid
when $F$ induces an isomorphism between the spaces of
invariant probability measures
(see Definition 3.1 and Corollary 3.11 of \cite{LM2}).

We denote the set of isometric homeomorphisms on $\T$ by $\Isom(\T)$.
The group $\Isom(\T)$ consists of the reflection and rotations.
For $t\in\T$ we write the translation $s\mapsto s+t$ on $T$ by $R_t$.
When $\xi:X\to\T$ is a continuous map,
$X\ni x\mapsto R_{\xi(x)}$ is a cocycle.
We denote this by $R_\xi$.

Let $\alpha\times\phi$ be a minimal homeomorphism on $X\times\T$ and
let $A=C^*(X\times\T,\alpha\times\phi)$.
Then $A$ is a unital simple $C^*$-algebra.
We will use $j_{\alpha}:C(X\times \T)\to A$ for the embedding
whenever it is convenient.

Let $u\in A$ be the implementing unitary.
For $x\in X$, let $A_x$ be the $C^*$-subalgebra generated by
$C(X\times\T)$ and $uC_0((X\setminus\{x\})\times\T)$.
By \cite[Proposition 3.3]{LM2},
$A_x$ is known to be a unital simple $A\T$ algebra and
the tracial state space $T(A_x)$ coincides with $T(A)$.
Besides, $A_x$ has real rank zero if and only if
$\alpha\times\phi$ is rigid.
We also remark that $A_x\cap C^*\xa$ is a unital simple AF algebra
(see \cite{Pu1}).

\section{Approximate unitary equivalence of unitaries}

The following is quoted from \cite{EGL2}.

\begin{lem}[{\cite[Theorem 2.11]{EGL2}}]\label{ILgong}
Let $F\subset C(\T)$ be a finite subset and $\ep>0.$
There exists $\eta_1>0$ with the property described as follows.

For any $\dt_1>0,$ there exist a positive integer $K$ and
a number $\eta_2>0$ such that
for any $\dt_2>0,$ there exist a finite subset
$H\subset C(\T)_{s.a.}$ and a positive integer
$N$ satisfying the following
condition.

if $\phi, \psi: C(\T)\to B,$ where
$B=\oplus_{j=1}^mM_{r(j)}(C([0,1]))$ are two unital \hm s such that
\begin{enumerate}
\item $\phi$ has the property $\sdp(\eta_1/32,\dt_1)$ and
$\sdp(\eta_2/32,\dt_2);$
\item $\lvert\tau\circ\phi(f)-\tau\circ\psi(f)\rvert<\dt_2/4$
for all $f\in H$ and all $\tau\in T(B);$
\item $r(j)\ge N,$ $j=1,2,\dots,m;$
\item $D(\phi(z)\psi(z)^*)\le 1/8K,$
\end{enumerate}
then there exists a unitary $w\in B$ such that
$$
\lVert\phi(f)-u^*\psi(f)u\rVert<\ep
$$
for all $f\in F.$
\end{lem}

\begin{rem}\label{IRgong}
Note that $K_1(B)=\{0\}.$
So the condition (5) in Theorem 2.11 of \cite{EGL2}
that $\phi_{*1}=\psi_{*1}$ is not needed.
In Theorem 2.11 of \cite{EGL2}, in this case, $B=M_r(C([0,1])$ is
a single summand and (3) should be replaced by $r\ge N.$
It is obvious that it works for finitely many summands
as long as $r(j)\ge N$ for all $j=1,2,\dots,m.$

We actually use a very special case of Theorem 2.11 of \cite{EGL2}.
A shorter proof could be given here.
In fact a version of this could be found in \cite{NT}.
We quote Theorem 2.11 of \cite{EGL2} for the convenience.
\end{rem}

Let $A$ be a unital $C^*$-algebra and
let $\phi:C(\T)\to A$ be a unital completely positive map.
For a subset ${\cal F}\subset C(\T)\setminus\{0\}$ and
a map $T=N\times K:{\cal F}\to\N\times\R_+$,
we say that $\phi$ is $T$-${\cal F}$-full,
if the following holds: for every $a\in{\cal F}$ there exist
$x_1,x_2,\dots,x_{N(a)}\in A$ such that
\[ \sum_{i=1}^{N(a)}x_i^*\phi(a)x_i=1 \]
and $\lVert x_i\rVert\le K(a)$ for all $i=1,2,\dots,N(a)$.
We say a unitary $u\in A$ is $T$-${\cal F}$-full,
if the homomorphism $f\mapsto f(u)$ is $T$-${\cal F}$-full.

\begin{lem}\label{ILsdp}
Let $A$ be a unital simple \CA\, and let $u\in A$ be a unitary.
Let $\eta>0.$
Suppose that $\{\zeta_1,\zeta_2,...,\zeta_m\}$ is an $\eta/128$-dense
subset of $\T,$ and that $g_i\in C(\T)$ satisfies
$0\le g_i\le 1,$ $g_i(t)=1$ if $t\in O(\zeta_i,\eta/512)$ and
$g_i(t)=0$ if $t\in\T\setminus O(\zeta_i,\eta/256).$
Suppose also that there are $x_{i,j}\in A$ such that
$$
\lVert\sum_{j=1}^{m(i)}x_{i,j}^*g_i(u)x_{i,j}
-1_A\rVert<1/4, \quad i=1,2,\dots,m.
$$
Then, for any $\ep>0$,
there are a finite subset ${\cal G}\subset A$ and $\gamma>0$
satisfying the following:

if $L:A\to B=\oplus_{k=1}^KM_{r(k)}(C([0,1]))$ is a unital
${\cal G}$-$\gamma$-multiplicative \morp\,
then there exists $\phi:C(\T)\to B$ such that
$$
\lVert\phi(z)-L(u)\rVert<\ep
$$
and $\phi$ has $\sdp(\eta/64,\dt)$ property,
where $\dt=1/\max\{m(i):1\le i\le m\}.$
\end{lem}

\begin{proof}
Since $C(\T)$ is semiprojective,
for the finite subset ${\cal F}=\{g_i: i=1,2,...,m\}\cup\{z\}$
and $\ep>0,$
there exist a finite subset ${\cal G}_1\subset C(\T)$ and
$\gamma>0$ such that,
for any unital ${\cal G}_1$-$\gamma$-multiplicative \morp\,
$L:C(\T)\to B,$ (where $B$ is any \CA\,)
there exists a unital \hm\, $\phi:C(\T)\to B$ for which
$$
\|\phi(f)-L(f)\|<\ep/2
$$
for all $f\in {\cal F}.$
Therefore, with a sufficiently large finite subset
${\cal G}$ containing $x_{i,j}$'s, $g_i(u)$'s and ${\cal G}_1$
and sufficiently small $\gamma,$ one has
$$
\|\sum_{j=1}^{m(i)} L(x_{i,j})^*\phi(g_i)L(x_{i,j})-1_B\|<1/2.
$$
For each $x\in\T,$ there is $\zeta_i$ such that
$O(x,\eta/64)\supset O(\zeta_i,\eta/128).$
It is then easy to check that
$\phi$ has the property $\sdp(\eta/32,\dt),$
where $\dt=1/\max\{m(i): i=1,2,\dots,m\}.$
\end{proof}

\begin{lem}\label{ILtrace}
Let $A$ be a unital stably finite simple \CA\,
and let ${\cal F}_1, {\cal F}_2\subset A_+$ be finite subsets.
Let $1/2>\gamma_1>0.$
Suppose that there exists a map $\pi:{\cal F}_1\to{\cal F}_2$
such that
$$
\lvert\tau(a)-\tau(\pi(a))\rvert<\gamma_1/8
$$
for all $\tau\in T(A).$
Then there exist $\dt>0$ and a finite subset ${\cal G}\subset A$
satisfying the following:

if $L:A\to B,$ where $B$ is a unital \CA\, with stable rank one,
is a unital ${\cal G}$-$\dt$-multiplicative \morp, then
$$
\lvert\tau'(L(a))-\tau'(L(\pi(a)))\rvert<\gamma_1
$$
for each $a\in {\cal F}_1$ and all $\tau'\in T(B).$
\end{lem}

\begin{proof}
It follows from \cite{CP} that, for each $a\in {\cal F}_1,$
there are $x_1(a),x_2(a),\dots,x_{n(a)}(a)\in A$ such that
$$
\|\sum_{k=1}^{n(a)}(x_k(a))^*x_k(a)-a\|<\gamma_1/4
\text{ and }\|\sum_{k=1}^{n(a)}x_k(a)(x_k(a))^*-\pi(a)\|<\gamma_1/4.
$$
Thus, with sufficiently large ${\cal G}$ and
sufficiently small $\dt>0,$ one has
\[ \|\sum_{k=1}^{n(a)}L(x_k(a))^*L(x_k(a))-L(a)\|<\gamma_1/2 \]
and
\[ \|\sum_{k=1}^{n(a)}L( x_k(a))L(x_k(a))^*-L(\pi(a))\|<\gamma_1/2. \]
It follows that
$$
|\tau(L(a))-\tau(L(b(a))|<\gamma_1
$$
for each $a\in {\cal F}_1$ and all $\tau\in T(B).$
\end{proof}

\begin{lem}\label{Ilcu}
Let $A$ be a unital simple \CA\, and let $u,v\in A$ be two unitaries.
Suppose that
$$
\dist({\bar u},{\bar v})<d\text{ in }U(A)/CU(A).
$$
Then, there exist $\dt_1>0,$  $\dt_2>0$ and a finite subset
${\cal G}\subset A$ satisfying the following:

If $L:A\to B=\sum_{j=1}^nM_{r(j)}(C([0,1]))$
(for any integer $n$ and $r(j)>0$)
is a unital ${\cal G}$-$\dt_2$-multiplicative \morp\, then
$$
D(u_1^*v_1)<2d,
$$
where $u_1,v_1$ are any unitaries in $B$ for which
$$
\|u_1-L(u)\|<\dt_1\text{ and }\|v_1-L(v)\|<\dt_1.
$$
\end{lem}

\begin{proof}
There are unitaries $a_1,a_2,\dots,a_m,b_1,b_2,\dots,b_m\in A$
such that
$$
\|u^*v-c\|<d,
$$
where $c=\prod_{j=1}^ma_jb_ja_j^*b_j^*.$
Choose $\dt_1=d/4.$
For sufficiently large ${\cal G}$ and
sufficiently small $\dt_2>0,$ one has
$$
\|L(u)^*L(v)-L(c)\|<d+d/4
$$
and there are unitaries
$a_1',a_2',\dots,a_m',b_1',b_2',\dots,b_m'\in B$ such that
$$
\|L(c)-\prod_{j=1}^ma_j'b_j'(a_j')^{-1}(b_j')^{-1}\|<d/4.
$$
If $\|u_1-L(u)\|<\dt_1$ and $\|v_1-L(v)\|<\dt_1$, then
$$
\|u_1^*v_1-c'\|<2d,
$$
where $c'=\prod_{j=1}^ma_j'b_j'(a_j')^{-1}(b_j')^{-1}.$
Clearly in $B,$ $D(c')=0.$
It follows that
\[ D(u_1^*v_1)<2d. \]
\end{proof}

\begin{lem}\label{ILuapp}
For any $\ep>0,$ $l\ge 2\pi$ and
$T:C(\T)_+\setminus\{0\}\to\N\times\R_+,$
there exist a finite subset
${\cal G}\subset C(\T)_+\setminus\{0\}$ and an integer $L>0$
satisfying the following:

Let $A$ be a unital simple \CA\, and
let $u,v\in A$ be unitaries such that
$[u]=[v]$ in $U(A)/U_0(A)$ and $\cel(u^*v)\le l$.
For any \hm\, $\phi: C(\T)\to A$ which is $T$-${\cal G}$-full,
there is a unitary $w\in M_{L+1}(A)$ such that
$$
\|w^*\diag(u,u_0,\dots,u_0)w-\diag(v,u_0,\dots,u_0)\|<\ep.
$$
where $u_0=\phi(z).$
\end{lem}

\begin{proof}
Note that $K_i(\T)=\Z,$ $i=0,1,$ and
they are generated by $1_{C(\T)}$ and $z$
(the identity map on $\T$).

Suppose that the lemma fails.
Then there would be an $\ep_0>0,$ $l_0 \ge 2\pi$ and
$T: C(\T)_+\setminus\{0\}\to\N\times\R_+$
such that the assertion does not hold.
Let $\{O_n\}_{n\in\N}$ be an open base of $\T$ and
let $g_n\in C(\T)_+$ be a function
satisfying $\{t:g_n(t)>0\}=O_n$.
Put ${\cal G}_n=\{g_1,g_2,\dots,g_n\}$.
Then we would have
two sequences of unitaries $\{u_n\},\{v_n\}$
in a sequence of unital simple \CA s $A_n$ for which $[u_n]=[v_n]$
in $U(A_n)/U_0(A_n)$ and $\cel(u_n^*v_n)\le l_0$, and
a $T$-${\cal G}_n$-full homomorphism $\phi_n:C(\T)\to A_n$
such that
$$
\inf\{\|w^*\diag(u_n,\zeta_0,\dots,\zeta_0)w-
\diag(v_n,\zeta_0,\dots,\zeta_0)\|:w\in M_n(A_n)\}\ge \ep_0,
$$
where $\zeta_0=\phi_n(z)$.
Let $B=\ell^{\infty}(\{A_n\}),$ $B_0=c_0(\{A_n\}),$
$U=\{u_n\}$, $V=\{v_n\}\in B,$
$\Phi_0: C(\T)\to B$ be defined by $\Phi_0(f)=\{\phi_n(f)\}$
for $f\in C(\T),$ and $\pi:B\to B/B_0$ be the quotient map.
Since
$$
\cel(u_n^*v_n)\le l_0,
$$
it is easy to see that
$$
\cel(U^*V)\le l_0\text{ in }B\text{ and }
\cel(\pi(U)^*\pi(V))\le l_0\text{ in }B/B_0.
$$
This, in particular, implies that $[\pi(U)]=[\pi(V)].$
Since $\phi_n$ is $T$-${\cal G}_n$-full,
we can see that $\pi\circ \Phi_0$ is full.
Thus, by Theorem 1.2 in \cite{GL}, there exists an integer $N>0$
and a unitary $W\in M_{N+1}(B/B_0)$ such that
$$
\|W^*\diag(U,\pi\circ\Phi_0(z),\dots,\pi\circ\Phi_0(z))W
-\diag(V,\Phi_n(z),\dots,\Phi_0(z))\|<\ep_0/2.
$$
Note that there exists a sequence of unitaries
$w_n\in M_{N+1}(A_n)$ such that $\pi(\{w_n\})=W.$
It follows that
$$
\|w_n^*\diag(u_n,\phi_n(z),\dots,\phi_n(z))w_n
-\diag(v_n,\phi_n(z),\dots,\phi_n(z))\|<\ep_0/2
$$
for all sufficiently large $n.$
This contradicts the assumption that the lemma fails.
\end{proof}

\begin{lem}\label{ILfull}
Let $A$ be a unital simple \CA\, and
let $\phi: C(\T)\to A$ be a monomorphism.
Suppose that $\phi$ is $T$-${\cal F}$-full,
where ${\cal F}$ is a finite subset of $C(\T)_+\setminus\{0\}$
and $T=N\times K:{\cal F}\to\N\times\R_+$ is a map.
Put $T'(f)=(N(f),2K(f))$ for $f\in{\cal F}$.
Then there exist $\dt>0$ and a finite subset ${\cal G}\subset A$
such that if $L:A\to B$ (for any unital \CA\,) is
a unital ${\cal G}$-$\dt$-multiplicative \morp\, and
$u\in B$ is a unitary satisfying $\|u-L\circ\phi(z)\|<\dt$,
then $u$ is $T'$-${\cal F}$-full.
\end{lem}

\begin{proof}
By assumption, for each $f\in{\cal F},$ there are
$x_1(f),x_2(f),\dots,x_{N(f)}(f)\in A$ such that
$$
\sum_{i=1}^{N(f)}(x_i(f))^*\phi(f)x_i(f)=1_A
$$
and $\|x_i(f)\|\le K(f)$.

It is clear that, with a sufficiently large ${\cal G}$ and
sufficiently small $\dt>0,$
$$
\|\sum_{i=1}^{N(f)}L(x_i(f))^*L\circ\phi(f)L(x_i(f))-1_B\|<1/4,
$$
provided that $L$ is a unital ${\cal G}$-$\dt$-multiplicative
\morp.
Then there exists $b(f)\in B_+$ with $\|b(f)\|<4/3$ such that
$$
\sum_{i=1}^{N(f)}b(f)L(x_i(f))^*L\circ\phi(f)L(x_i(f))b(f)=1_B.
$$
Note that $\|L(x_i(f))b(f)\|\le 4K(f)/3.$
By choosing a small $\dt$, we may assume that
$L\circ\phi(f)$ is sufficiently close to $f(u)$
for every $f\in{\cal F}$, and
$$
\|\sum_{i=1}^{N(f)}b(f)L(x_i(f))^*f(u)L(x_i(f))b(f)-1_B\|<1/4.
$$
By repeating the same argument as above we can conclude that
$u$ is $T'$-${\cal F}$-full.
\end{proof}

\begin{lem}\label{ILdiv}
Let $A$ be a unital simple \CA\, with $TR(A)\le 1.$
Then, for any $\ep>0,$ any $\sigma>0,$ any integers $m_0, N>0$ and
any finite subset ${\cal F}\subset A,$ there exist mutually
orthogonal projections $q,p_1,p_2,\dots,p_m$ ($m\ge m_1$) with
$[q]\le [p_1]$ and $[p_1]=[p_i],$ $i=1,2,\dots,m,$
a \SCA\, $C\in {\cal I}$ with $1_C=p_1$ for which
each summand of $C$ has rank at least $N,$ and
unital ${\cal F}$-$\ep$-multiplicative \morp s
$L_1:A\to qAq$ and $L_2:A\to C$ such that
$$
\|x-L_1(x)\oplus\diag(L_2(x),L_2(x),\dots,L_2(x))\|<\ep,
$$
where $L_2$ is repeated $m$ times, for all $x\in {\cal F},$
$\tau(q)<\sigma$ and $2\tau(q)>\tau(p_1)$ for all $\tau\in T(A).$
\end{lem}

\begin{proof}
The proof is a minor modification of that of Lemma 5.5
in \cite{Lntain}.
The only difference is that, in Lemma 5.5 of \cite{Lntain},
one does not have $2\tau(q)>\tau(p_1)$ (for all $\tau\in T(A)$).

Choose $m\ge m_1$ so that $1/m<\sigma/2.$
In the proof of Lemma 5.5 in \cite{Lntain}, choose $n=2m+1.$
Then, the statement of the lemma holds,
where $L_2$ is repeated $2m+1$ times and
$\tau(q)<\sigma/2$ (for all $\tau\in T(A)$).
Then one replaces $L_1$ by $L_1\oplus L_2$ and
$L_2$ by $L_2\oplus L_2.$
Then the present lemma holds.
\end{proof}

\begin{lem}\label{ILcut}
Let $A$ be a unital simple \CA\, with property (SP) and
let $u, v\in A$ be two unitaries with $1\in sp(u), sp(v).$
Then, for any $\ep>0,$ there is a unitary $w\in A$ and
a nonzero projection $e\in A$
and two unitaries $u_1,v_1\in (1-e)A(1-e)$ such that
$$
\|(e+u_1)-u\|<\ep/2\text{ and }\|(e+v_1)-w^*vw\|<\ep/2.
$$
\end{lem}

\begin{proof}
Let $\dt>0.$ Define $f\in C(\T,[0,1])$ such that
$f(\xi)=1$ if $\dist(\xi,1)<\dt/2$ and
$f(\xi)=0$ if $\dist(\xi,1)\ge \dt.$
Since $A$ has property (SP), there are nonzero projections
$e_1\in \overline{f(u)Af(u)}$ and $e_2\in \overline{f(v)Af(v)}.$
Since $A$ is simple there exist nonzero projections
$e_i'\le e_i$ such that $e_1'$ is unitarily equivalent to $e_2'$.
Put $e=e_1'.$
If $\dt$ is sufficiently small, we have
$$
\|eu-ue\|<\ep/8\text{ and }\|e_2'v-ve_2'\|<\ep/8.
$$
It is easy to obtain unitaries $u_1\in (1-e)A(1-e)$ and $u_2\in
(1-e_2')A(1-e_2')$ such that
$$
\|(e+u_1)-u\|<\ep/2\text{ and }\|(e_2'+u_2)-v\|<\ep/2.
$$
There is a unitary $w\in A$ such that $w^*e_2'w=e.$
Then, with $v_1=w^*u_2w,$
$$
\|(e+v_1)-w^*vw\|<\ep/2.
$$
\end{proof}

\begin{lem}\label{ILcel}
Let $A$ be a unital simple \CA\, and $u,v\in U(A).$
Suppose that $\cel(u^*v)\le l$ for some $l>0.$
Then, for any $\ep>0,$ there is a finite subset ${\cal F}\subset A$
and $\dt>0$ such that,
for any unital ${\cal F}$-$\dt$-multiplicative \morp\, $L:A\to B$
(for any unital \CA\,$B$),
there are two unitaries $u_1,v_1\in B$ such that
$$
\|L(u)-u_1\|<\ep/2,\quad\|L(v)-v_1\|<\ep/2\text{ and}
$$
$$
\cel(u_1^*v_1)<l+\ep.
$$
\end{lem}

\begin{proof}
This is essentially the same statement of Lemma 6.8 of \cite{Lntain}.
\end{proof}

Now we are ready to prove the following theorem.

\begin{thm}\label{Tu1}
Let $\ep>0$ and let $T:C(\T)_+\setminus\{0\}\to\N\times\R_+$ be a map.
Then there exist $\dt>0$ and
a finite subset ${\cal F}\subset C(\T)_+\setminus\{0\}$
satisfying the following:
For any unital simple \CA\, with $TR(A)\le 1$ and
any $T$-${\cal F}$-full unitary $u\in U(A)$ with $sp(u)=\T$,
if $v\in U(A)$ is a unitary such that
$$
[u]=[v]\text{ in }K_1(A),\quad\dist({\bar u},{\bar v})<\dt
$$
and
$$
|\tau\circ f(u)-\tau\circ f(v)|<\dt
$$
for all $f\in{\cal F}$ and all $\tau\in T(A),$
then there is a unitary $w\in A$ such that
$$
\|w^*uw-v\|<\ep.
$$
\end{thm}

\begin{rem}
In the statement above, for any $u\in A$ with $sp(u)=\T,$
since $A$ is a simple unital \CA\, the map $T$ always exists.
We would like to point out that $\dt$ and ${\cal F}$
depend only on such $T$ but not on $A$ or the choice of $u$
as long as the map $T$ works for $u.$
\end{rem}

\begin{proof}
Let $l=9\pi.$
Suppose that $T:C(\T)_+\setminus\{0\}\to\N\times\R_+$ is defined
by $T(f)=(N(f), K(f)).$
By applying Lemma \ref{ILuapp} with $\ep/16,$ $l$ and
$T'(f)=(N(f),2K(f)),$
we get an integer $L>0$ and
a finite subset ${\cal G}\subset C(\T)_+\setminus\{0\}.$

Let $\eta_1>0$ be as in Lemma \ref{ILgong} corresponding to $\ep/4$
and $\{z\}\subset C(\T)$.
Let $\{\zeta_1,\zeta_2,\dots,\zeta_m\}\subset\T$ be
an $\eta_1/128$-dense subset of $\T.$
Choose functions $g_i\in C(\T)$ such that $0\le g(t)\le 1,$
$g_i(t)=1$ if $t\in O(\zeta_i,\eta_1/512)$ and $g_i(t)=0$
if $t\in \T\setminus O(\zeta_i,\eta_1/256),$ $i=1,2,\dots,m.$
Let $\dt_1=1/\max\{N(g_i):1\le i\le m\}.$
By using Lemma \ref{ILgong} with $\{z\},$ $\ep/4,$ $\eta_1$ and
$\dt_1,$
we obtain a natural number $K$ and $\eta_2>0$.

Now let $\{\xi_1,\xi_2,\dots,\xi_l\}\subset\T$ be
$\eta_2/128$-dense in $\T.$
Let $h_i\in C(\T)$ such that $0\le h_i(t)\le 1,$
$h_i(t)=1$ if $t\in O(\xi_i,\eta_2/512)$ and $h_i(t)=0$
if $t\in \T\setminus O(\xi_i,\eta_2/256),$ $i=1,2,\dots,l.$
Let $\dt_2=1/\max\{N(h_i):1\le i\le l\}.$
Let $H$ and $N>0$ be described in Lemma \ref{ILgong}
corresponding to the above
$\ep/4,$ $\eta_1>0,$ $\dt_1,$ $K,$ $\eta_2$ and $\dt_2.$
Set ${\cal H}_1=\{a+, a_-,a: a\in H\}.$

Now let ${\cal F}$ be a finite subset which contains
${\cal G},$ $H_1$ above and $\{g_1,\dots,g_m,h_1,h_2,\dots,h_l\}.$
We may assume that $K>16.$
So we choose $\dt=\min\{\pi/16K,\dt_2/64\}.$
Note that ${\cal F}$ and $\dt$ depend only on $T$ and $\ep$.

We would like to show that ${\cal F}$ and $\dt$ do the work.
Suppose that $A$ is a unital simple $C^*$-algebra with $TR(A)\le 1$
and $u\in U(A)$ is $T$-${\cal F}$-full and $sp(u)=\T$.
Let $v\in A$ be another unitary such that $[u]=[v]$ in $K_1(A)$,
$$
\dist({\bar u},{\bar v})<\dt\le\pi/16K
$$
and
$$
|\tau\circ f(u)-\tau\circ f(v)|<\dt\le\dt_2/64
$$
for all $f\in H_1\subset{\cal F}$ and all $\tau\in T(A).$

By applying Lemma \ref{ILcut}, we may assume, without loss of generality,
$u=e+u'$ and $v=e+v',$ where $u',v'\in(1-e)A(1-e)$ are two unitaries.
To simply notation, fix a nonzero projection $e\in A$ so that
it suffices to prove the following:
there exists a unitary $w\in (e+1_A)M_2(A)(e+1_A)$ such that
$$
\|w^*(e+u)w-(e+v)\|<\ep.
$$

It follows from Lemma 6.9 in \cite{Lntain} that
$$
\cel(u^*v)<8\pi+\pi/16K<9\pi=l.
$$
Define $\phi,\psi:C(\T)\to A$ by $\phi(f)=f(u)$ and $\psi(f)=f(v)$
for all $f\in C(\T).$

Since ${\cal F}$ contains $\{g_1,g_2,\dots,g_m\}$,
we have
$$
\sum_{j=1}^{N(g_i)}(x_j(g_i))^*\phi(g_i)x_j(g_i)=1_A,
$$
for some $x_j(g_i)\in A$ with
$\|x_j(g_i)\|\le K(g_i),$ $i=1,2,\dots,m.$
Since ${\cal F}$ contains $\{h_1,h_2,\dots,h_l\}$,
we also have
$$
\sum_{j=1}^{N(h_i)}(x_j(h_i))^*\phi(h_i)x_j(h_i))=1_A
$$
for some $x_j(h_i)\in A$ with $\|x_j(h_i)\|\le K(h_i),$
$i=1,2,\dots,l.$

Now we apply Lemma \ref{ILdiv}.
For a finite subset ${\cal G}_0\subset A$ and $\dt_0>0,$
we have a natural number $L'$ greater than $L$ and
$$
\|x-(L_1(x)\oplus
\diag(\overbrace{L_2(x),L_2(x),\dots,L_2(x)}^{L'})\|<\ep/16
$$
for both $x=u,v,$ where $L_1:A\to qAq$ and $L_2:A\to B$ are
${\cal G}_0$-$\dt_0$-multiplicative \morp s,
where $B=\oplus_{j=1}^mM_{r(j)}(C([0,1]))\subset p_1Ap_1$
with $r(j)\ge N$ ($1\le j\le m$), $1_B=p_1,$
$[q]\le [p_1]$, $2[q]\ge [p_1]$ and $[q] \le [e]$.
We choose ${\cal G}_0$ so large and $\dt_0$ so small that,
by Lemma \ref{ILsdp}, there is a \hm\, $\phi_1:C(\T)\to B$ so that
$$
\|\phi_1(f)-L_2(\phi(f))\|<\min\{\ep/4,\dt_2/16\}
$$
for all $f\in H_1\cup\{z\}$
and $\phi_1$ has $\sdp(\eta_1/64,\dt_1)$
and $\sdp(\eta_2/64,\dt_2)$ property.
Moreover, since $\phi$ is $T$-${\cal G}$-full,
by Lemma \ref{ILfull}, we may assume that
$\phi_1$ is $T'$-${\cal G}$-full.

We may also assume that
there is a \hm\, $\psi_1:C(\T)\to B$ such that
$$
\|\psi_1(f)-L_2(\psi(f))\|<\min\{\ep/4,\dt_2/16\}
$$
for all $f\in H_1\cup\{z\}$.

By Lemma \ref{Ilcu}, we may assume that
$$
D(\phi_1(z)^*\psi_1(z))<1/8K.
$$

With sufficiently large ${\cal G}_0$ and sufficiently small $\dt_0,$
by Lemma \ref{ILtrace}, we may assume that
$$
|\tau'(\psi_1(f))-\tau'(\phi_1(f))|<\dt_2/4
$$
for all $f\in H$ and all $\tau'\in T(B).$
Now by applying Lemma \ref{ILgong},
we obtain a unitary $w_1\in B$ such that
$$
\|w_1^*\phi_1(z)w_1-\psi_1(z)\|<\ep/4.
$$

In the above, for sufficiently large ${\cal G}_0$ and
sufficiently small $\dt_0,$ by Lemma \ref{ILcel}, we may also assume
that there are unitaries $u_2, v_2\in qAq$ such that
$$
\|L_1(u)-u_2\|<\ep/16,\quad\|L_1(v)-v_2\|<\ep/16
$$
and
$$
\cel(u_2^*v_2)\le 9\pi.
$$

Since $[q]\le[p_1]$ and $[p_1]-[q]\le[q]\le[e]$,
there is $e_1\le e$ such that $[e_1]+[q]=[p_1].$
Let $u_3=e_1+u_2$ and $v_3=e_1+v_2.$

Put
$$
U=\diag(\overbrace{\phi_1(z),\phi_1(z),\dots,\phi_1(z)}^{L'})
$$
and
$$
V=\diag(\overbrace{\psi_1(z),\psi_1(z),\dots,\psi_1(z)}^{L'}).
$$
Then by the choice of $L,$ by applying Lemma \ref{ILuapp}, we obtain
a unitary $w_2\in (e_1+1_A)M_2(A)(e_1+1_A)$ such that
$$
\|w_2^*(u_3\oplus U)w_2-(v_3\oplus U)\|<\ep/16.
$$
Let $w_3=w_2(q+\diag(\overbrace{w_1,w_1,\dots,w_1}^{L'})).$
It follows that
$$
\|w_3^*(u_3\oplus U)w_3-(v_3\oplus V)\|<\ep/16+\ep/4.
$$
Combining all the above, we have
$$
\|w_3^*(e_1+u)w_3-(e_1+v)\|<3(\ep/16+\ep/4)<\ep.
$$
Set $w=(e-e_1)+w_3.$ Then
$$
\|w^*(e+u)w-(e+v)\|<\ep.
$$
\end{proof}

\begin{cor}\label{1C}
Let $A$ be a unital simple \CA\, with $TR(A)\le 1$ and
let $u,v\in A$ be two unitaries with $sp(u)=sp(v)=\T.$
Then there exists a sequence of
unitaries $w_n\in A$ such that
$$
\lim_{n\to\infty}w^*_nuw_n=v
$$
if and only if
$$
[u]=[v]\text{ in }K_1(A),\quad{\bar u}={\bar v}
\text{ in }U(A)/CU(A)
$$
and
$$
\tau(f(u))=\tau(f(v))
$$
for all $f\in C(\T)$ and all $\tau\in T(A).$
\end{cor}

\begin{lem}\label{IIILetoI}
Let $e$ be a nonzero projection of a unital simple \CA\, $A$ with
$TR(A)\le 1.$
Then $\imath: U(eAe)/CU(eAe)\to U(A)/CU(A)$ defined
by ${\bar v}\to \overline{v+(1-e)}$ is a continuous isomorphism.
\end{lem}

\begin{proof}
It is clear that  $\imath$ is a continuous \hm.
It follows from  Theorem 6.7 in \cite{Lntain} that it is surjective.
Suppose that ${\bar u}\in {\rm ker}\imath.$ Thus $u+(1-e)\in CU(A).$
It follows Lemma 6.9 in \cite{Lntain} that $u +(1-e)\in U_0(A).$
Since $A$ has stable rank one, by \cite{Rf},
it is easy to see that $[u+(1-e)]=0$ in $K_1(A).$
Since $A$ is simple, we conclude that
$[u]=0$ in $K_1(eAe).$ Since $A$ has stable rank one,
it follows that $u\in U_0(eAe).$
By expressing $u$ as finite product of exponentials, we obtain
a piecewise smooth map $\eta: [0,1]\to U_0(eAe)$ with
$\eta(0)=e$ and $\eta(1)=u.$
Define $\xi:[0,1]\to U(A)$ by $\xi(t)=\eta(t)+(1-e).$
Then
$$
\delta_A(\xi)(\tau)=\frac{1}{2\pi i}\int_0^1\tau(\xi'(t)\xi(t)^*)\,dt
=\frac{1}{2\pi i}\int_0^1\tau(\eta'(t)\eta(t)^*)\,dt
$$
for all $\tau\in T(A).$ The fact that $u+(1-e)\in CU(A)$ implies
that $\delta_A(\xi)\in \overline{\rho_A(K_0(A))}$ (see \cite{Th}).
Suppose that there are $x_n\in K_0(A)$ such that
$\tau(x_n)\to \delta(\xi)(\tau)$ uniformly on $T(A).$
Then
$$
\tau(x_n)/\tau(e)\to \delta_A(\xi)(\tau)/\tau(e).
$$
For each $\tau\in A,$ define $\tilde\tau(a)=\tau(a)/\tau(e)$
for $a\in eAe.$
So $\delta_{eAe}(\eta)(\tilde\tau)=\delta_A(\xi)(\tau)/\tau(e).$
Since $K_0(eAe)=K_0(A),$ we conclude that
$\delta_{eAe}(\eta)\in \overline{\rho_A(K_0(eAe))}.$
Equivalently $u\in CU(eAe).$
Thus $\imath$ is injective.
\end{proof}

\begin{thm}\label{ITX}
Let $A$ be a unital simple \CA\, with $TR(A)\le 1$ and
$X$ be the Cantor set.
Then two unital monomorphisms $h_1,h_2:C(X\times\T)\to A$
are approximately unitarily equivalent if and only if
\[ (h_1)_{*i}=(h_2)_{*i}, \ i=0,1,\quad
h_1^{\sharp}=h_2^{\sharp} \]
and
\[ \tau\circ h_1(f)=\tau\circ h_2(f) \]
for all $f\in C(X\times\T)$ and $\tau\in T(A).$
\end{thm}

\begin{proof}
The ``only if" is clear. We will show the ``if" part.

Let $\ep>0$ and ${\cal F}\subset C(X\times \T)$ be a finite subset.
Without loss of generality, we may assume that
\[ {\cal F}=\{f_i,f_i\times z: i=1,2,\dots,m\}, \]
where $f_i=1_{O_i}$ and $O_1,O_2,\dots,O_m$ are mutually disjoint clopen
subset of $X$ for which $\cup_{i=1}^mO_i=X.$

Since $(h_1)_{*i}=(h_2)_{*i}$, $i=0,1,$
there is a unitary $w_1\in A$ such that
$$
w^*h_1(f_i)w=h_2(f_i),\quad i=1,2,\dots,m.
$$
To simplify notation, we may assume that
$h_1(f_i)=h_2(f_i)=p_i,$ $i=1,2,\dots,m.$
Working in each $p_iAp_i,$ by applying Lemma \ref{IIILetoI} and
Corollary \ref{1C}, there is a unitary $u_i\in p_iAp_i$ such that
$$
\|u_i^*h_1(f_i\times z)u_i-h_2(f_i\times z)\|<\ep,
$$
for all $i=1,2,\dots,m.$
Define $w_2=\sum_{i=1}^mu_i.$ Then $w_2\in U(A)$ and
$$
\|w_2^*h_1(g)w_2-h_2(g)\|<\ep
$$
for all $g\in{\cal F}.$
\end{proof}

Combing the proof of Theorem \ref{ITX} and \ref{Tu1},
we actually prove the following.

\begin{cor}\label{ITX2}
Let $A$ be a unital simple \CA\, with $TR(A)\le 1$ and
$X$ be the Cantor set.
Fix a monomorphism $h_1:C(X\times\T)\to A.$
Then, for any $\ep>0$ and
any finite subset ${\cal F}\subset C(X\times\T),$
there exist $\dt>0,$ a finite subset ${\cal G}\subset C(X\times\T),$
a finitely generated subgroup $G_0\subset K_0(C(X\times \T))$ and
a finitely generated subgroup $G_1\subset K_1(C(X\times \T))$
satisfying the following:
if $h_2:C(X\times\T)\to A$ is a monomorphism such that
$$
(h_2)_{*i}|_{G_i}=(h_1)_{*i}|_{G_i}, \quad
|\tau(h_2(g))-\tau(h_1(g))|<\dt
$$
for all $g\in {\cal G}$ and $\tau\in T(A)$, and
$$
\dist(h_1^{\sharp}(g),h_2^{\sharp}(g))<\dt
$$
for all $g\in G_1,$
then there exists a unitary $W\in A$ such that
$$
\lVert Wh_2(f)W^*-h_1(f)\rVert<\ep
\text{ for all }f\in {\cal F}.
$$
\end{cor}

\section{Tracial rank}

Let $\xa$ be a Cantor minimal system and let $\xi\in C(X,\T)$.
In this section, we will only consider the case that
$\alpha\times R_\xi$ is minimal.
Put $A=C^*(X\times\T,\alpha\times R_\xi)$.
The purpose of this section is to show that the tracial rank of
$A$ is no more than one.
We denote the implementing unitary of $A$ by $u$.
In this section, we identify the circle $\T$
with the quotient space $\R/\Z$.

\begin{prop}\label{IIIP'}
Let $x\in X$ and let $U$ be a clopen neighborhood of $x\in X$.
Suppose that there exists $M\in\N$ such that
\[ \|u^Mzpu^{M*}-zq\|<\ep, \]
where $p=1_U$ and $q=u^Mpu^{M*}$.
Then there exists a partial isometry $w\in A_x$ such that
$w^*w=p$, $ww^*=q$ and
\[ \|wzpw^*-zq\|<\ep. \]
\end{prop}

\begin{proof}
There exists a unitary normalizer $w_1\in A_x\cap C^*\xa$ of $C(X)$
such that $w_1pw_1^*=q$.
We may assume that there exists a continuous function $n:X\to\Z$
such that $w_1=\sum_{k\in\Z}u^k1_{n^{-1}(k)}$.
Since $u^*zu=e^{2\pi\sqrt{-1}\xi}z$,
we can find a continuous map $\eta:U\to\T$ such that
\[ w_1^*u^Mzpu^{M*}w_1=e^{2\pi\sqrt{-1}\eta}zp. \]
Clearly we have $[w_1^*u^Mzpu^{M*}w_1]=[zp]$ in $K_1(pA_xp)$.
We also get $\tau(w_1^*u^Mzpu^{M*}w_1)=\tau(zp)$
for all $\tau\in T(pA_xp)$, because $T(A)\cong T(A_x)$.
Furthermore
\[ w_1^*u^Mzpu^{M*}w_1(zp)^*=e^{2\pi\sqrt{-1}\eta} \]
belongs to $B=p(A_x\cap C^*\xa)p$,
which is a unital simple infinite dimensional AF algebra.
Hence, the unitary $e^{2\pi\sqrt{-1}\eta}$ is contained
in $U(B)=CU(B)\subset CU(pA_xp)$.
Thus, Corollary \ref{1C} applies and yields a unitary
$w_2\in pA_xp$ such that
$$
\|w_1^*u^Mzpu^{M*}w_1-w_2zpw_2^*\|<\sigma,
$$
where $\sigma=\ep-\|u^Mzp^{M*}-zp\|.$
Put $w=w_1w_2$. Then
$$
\lVert wzpw^*-zq\rVert\le\lVert w_1(w_2zpw_2^*)w_1^*-u^Mzpu^{M*}\rVert
+\lVert u^Mzpu^{M*}-zq\rVert<\ep.
$$
\end{proof}

The following is an improvement of Lemma 5.5 of \cite{LM2}.
\begin{lem}\label{5L1}
Let $x\in X$.
For any $N\in\N$, $\ep>0$ and
a finite subset ${\cal F}\subset C(X\times\T)$,
we can find a natural number $M>N$,
a clopen neighborhood $U$ of $x$ and a partial isometry $w\in A_x$
which satisfy the following.
\begin{enumerate}
\item $\alpha^{-N+1}(U),\alpha^{-N+2}(U),\dots,
U,\alpha(U),\dots,\alpha^M(U)$
are mutually disjoint, and
$\mu(U)<\varepsilon/M$ for all $\alpha$-invariant measure $\mu$.
\item $w^*w=1_U$ and $ww^*=1_{\alpha^M(U)}$.
\item $u^{*i}wu^i\in A_x$ for all $i=0,1,\dots,N-1$.
\item $\lVert wf-fw\rVert<\varepsilon$ for all $f\in{\cal F}$.
\end{enumerate}
\end{lem}
\begin{proof}
Without loss of generality, we may assume
${\cal F}=\{f_1,f_2,\dots,f_k,z\}$,
where $f_i$ belongs to $C(X)\subset C(X\times\T)$.
There exists a clopen neighborhood $O$ of $x$ such that
\[ \lvert f_i(x)-f_i(y)\rvert<\varepsilon/2 \]
for all $y\in O$ and $i=1,2,\dots,k$. Since $\alpha\times R_\xi$
is minimal, we can find $M>N$ such that
$(\alpha\times R_\xi)^M(x,0)\in O\times I,$
where $I=\{t\in\T:|t|<\ep\}.$
Let $U$ be a clopen neighborhood of $x$ such that the condition (1) is
satisfied and
\[ (\alpha\times R_\xi)^M(y,0)\in O\times I \]
for all $y\in U$.
If $1/K<\ep,$
an easy way to get $\mu(U)<\ep/M$ is to choose $U$ so that\\
$\alpha^{-N+1}(U),\dots,U,\dots,\alpha^{M+K}(U)$ are mutually disjoint.
Moreover, we require that $U\cup\alpha^M(U)\subset O.$
Let $p=1_U$ and $q=1_{\alpha^M(u)}.$ Since
$(\alpha\times R_\xi)^M= \alpha^M\times R_{\eta}$ for some
$\eta\in C(X,\T),$ we check that
$$
\|u^Mzpu^{*M}-zq\|<\ep.
$$
By applying Lemma \ref{IIIP'}, we obtain a partial isometry $w\in A_x$
which satisfies (2) and
$$
\|wzpw^*-zq\|<\ep.
$$
Since $U\cup \alpha^M(U)\subset O,$ by the choice of $O,$
it is easy to check that
$$
\|wf_i-f_iw\|<\ep
$$
for all $i=1,2,\dots,k.$
To see (3), we note that
$$
pu^i=pu1_{\alpha^{-1}(U)}u1_{\alpha^{-2}(U)}\cdots
u1_{\alpha^{-i}(U)}
$$
and
$$
(u^{*i}q)^*=qu1_{\alpha^{M-1}(U)}u1_{\alpha^{M-2}(U)}\cdots
u1_{\alpha^{M-i}(U)}
$$
for every $i=1,2,\dots,N-1.$
Since $x\in U,$ by the condition (1),
one sees that $pu^i$ and $u^{*i}q$ belong to $A_x.$
It follows that $u^{*i}wu^i\in A_x$ for all $i=1,2,\dots,N-1.$
\end{proof}

\begin{thm}\label{5T1}
Let $\xa$ be a Cantor minimal system and
let $\xi:X\to\T$ be a continuous map.
If $\alpha\times R_\xi$ is minimal, then
$A=C^*(X\times\T,\alpha\times R_\xi)$ has tracial rank zero or one.
Consequently $A=C^*(X\times \T,\alpha\times R_{\xi})$ is
a unital simple $A\T$-algebra.
Moreover, it has tracial rank zero if and only if
$\alpha\times R_{\xi}$ is rigid.
\end{thm}

\begin{proof}
The proof is exactly the same as that of Theorem 5.6 of \cite{LM2}
when one uses Lemma \ref{5L1} instead of \cite[Lemma 5.5]{LM2}.
Only difference is that we do not assume that
$A_x$ has tracial rank zero. But we use the fact that
$A_x$ is a unital simple $A\T$-algebra
(see Proposition 3.3 of \cite{LM2}).

Let ${\cal F}\subset A$ be a finite subset and let $\ep>0.$
Fix $x\in X.$
By applying Lemma \ref{5L1},
exactly as in the proof of Theorem 5.6 of \cite{LM2},
one obtains a projection $e\in A_x$ such that the following hold.
\begin{itemize}
\item $\|ea-ae\|<\ep$ for all $a\in {\cal F}$.
\item For every $a\in {\cal F}$, there exists $b\in eA_xe$
such that $\lVert eae-b\rVert<\ep$.
\item $\tau(1-e)<\ep$ for all $\tau\in T(A).$
\end{itemize}
Since $A_x$ is a unital simple $A\T$-algebra (which has tracial
rank one or zero), using the fact that $A$ has stable rank one and
weakly unperforated $K_0(A)$ (see Theorem 3.12 in \cite{LM2})
and applying Theorem 4.8 in \cite{HLX},
exactly as in the proof of Theorem 5.6 in \cite{LM2},
we conclude that $A$ has tracial rank one or zero.

By Lemma 2.4 of \cite{LM2},
both $K_0(A)$ and $K_1(A)$ are torsion free.
It follows from \cite{Lntain} that
$A$ is isomorphic to a unital simple $A\T$-algebra.
\end{proof}

\section{Non-orientation preserving cases}

In this section we will show that the crossed product
$C^*(X\times\T,\alpha\times\phi)$ has tracial rank zero
if the cocycle $\phi$ takes its values in $\Isom(\T)$ and
$\alpha\times\phi$ is rigid.

\begin{lem}\label{rr0}
Let $A$ be a $C^*$-algebra with real rank zero and
let $E$ be a finite dimensional $C^*$-subalgebra with the same unit as $A$.
Then $B=A\cap E'$ also has real rank zero.
\end{lem}

\begin{proof}
Let $p_1,p_2,\dots,p_n$ be a family of mutually orthogonal central projections
of $E$ with $\sum p_i=1$.
Then $Ep_i$ is isomorphic to a full matrix algebra.
Since $p_i$ is central in $B$, it suffices to show that
$Bp_i$ has real rank zero for all $i=1,2,\dots,n$.
But this is obvious because $Bp_i=p_iAp_i\cap (Ep_i)'$ is isomorphic to
$e_iAe_i$ where $e_i$ is a minimal projection of $Ep_i$.
\end{proof}

Let $A$ be a unital $C^*$-algebra.
For $a\in A$, we define
\[ \lVert a\rVert_2=\sup_{\tau\in T(A)}\tau(a^*a)^{1/2}. \]
Then $\lVert\cdot\rVert_2$ is a norm on $A$.

\begin{lem}\label{CBS}
Let $e$ be a self-adjoint element of $A$ and
let $\{x_n\}_n$ be a sequence of self-adjoint elements of $A$.
Suppose that $\lim_{n\to\infty}\lVert x_n-e\rVert_2=0$ and
$\lVert e\rVert\leq1,\lVert x_n\rVert\leq1$ for all $n\in\N$.
Then, for every continuous function $f$ on $[-1,1]$, we have
\[ \lim_{n\to\infty}\lVert f(x_n)-f(e)\rVert_2=0. \]
\end{lem}

\begin{proof}
It suffices to show the claim when $f$ is a polynomial.
But this is obvious
because of $\lVert ab\rVert_2\leq\lVert a\rVert\lVert b\rVert_2$.
\end{proof}

\begin{lem}\label{replace}
Let $A$ be a unital simple $C^*$-algebra with tracial rank zero
and let $\{e_n\}_{n\in\N}$ be a sequence of projections in $A$
which satisfies
\[ \lim_{n\to\infty}\lVert ae_n-e_na\rVert_2=0 \]
for every $a\in A$.
Then there exist a subsequence $\{e_{m(n)}\}_{n\in\N}$ and
a sequence $\{x_n\}_{n\in\N}$ of projections in $A$
such that the following conditions are satisfied.
\begin{enumerate}
\item For every $a\in A$, we have $\lVert ax_n-x_na\rVert\to0$.
\item $\lVert e_{m(n)}-x_n\rVert_2\to0$.
\end{enumerate}
\end{lem}
\begin{proof}
Let $\{a_n\}_{n\in\N}$ be a dense sequence of $A$.
Since $A$ has tracial rank zero, we can find a projection $p_n$
and a unital finite dimensional $C^*$-algebra $E_n\subset p_nAp_n$
such that the following are satisfied.
\begin{itemize}
\item For every $i=1,2,\dots,n$,
$\rVert a_ip_n-p_na_i\lVert<1/n$.
\item For every $i=1,2,\dots,n$, there exists $b\in E_n$
such that $\rVert p_na_ip_n-b\lVert<1/n$.
\item $\lVert1-p_n\rVert_2<1/n$.
\end{itemize}
Using the Haar measure on the compact group $U(E_n),$
we define
\[ x_{m,n}=\int_{U(E_n)}ue_mu^*\,du. \]
Then $x_{m,n}\in A.$
It is then easy to check that $x_{m,n}$ commutes with
unitaries of $E_n$, and so it commutes with all elements of $E_n$.
Thus $x_{m,n}$ is a positive element lying in $p_nAp_n\cap E_n'$.
Hence, for every $i=1,2,\dots,n$, we have
$\lVert a_ix_{m,n}-x_{m,n}a_i\rVert<4/n$.
Moreover, by choosing a sufficiently large $m\in\N$,
we obtain
\[ \lVert e_m-x_{m,n}\rVert_2<\frac{1}{n}, \]
because $\lim_{m\to\infty}\lVert ue_m-e_mu\rVert_2=0$ for every $u\in U(E_n).$

In this way, we can find a subsequence $\{e_{m(n)}\}_n$ and
a sequence $\{x_n\}_n$ which satisfy the requirements (1) and (2).
It remains to replace $x_n$ to a projection.
Since $p_nAp_n\cap E_n'$ has real rank zero by Lemma \ref{rr0},
we may assume that $x_n$ has finite spectrum.
Let $f,g$ and $h$ be functions on $[0,1]$ defined by
\[ f(t)=\begin{cases} 0 & 0\leq t<1/2 \\
2t-1 & 1/2\leq t\leq 1, \end{cases}\]
\[ g(t)=\begin{cases} 2t & 0\leq t<1/2 \\
1 & 1/2\leq t\leq 1 \end{cases}\]
and $h(t)=1_{(1/2,1]}(t)$.
Then by using Lemma \ref{CBS}
we have $\lVert e_{m(n)}-f(x_n)\rVert_2\to0$ and
$\lVert e_{m(n)}-g(x_n)\rVert_2\to0$.
It follows from $0\leq(h-f)^2\leq(g-f)^2$ that
\[ \lim_{n\to\infty}\lVert h(x_n)-e_{m(n)}\rVert_2
=\lim_{n\to\infty}\lVert h(x_n)-f(x_n)\rVert_2
\leq\lim_{n\to\infty}\lVert g(x_n)-f(x_n)\rVert_2=0. \]
Since $h(x_n)$ still lies in $p_nAp_n\cap E_n'$,
it almost commutes with $a_1,a_2,\dots,a_n$.
Thus $h(x_n)$ is the desired projection.
\end{proof}

\begin{prop}\label{trRohlin}
Let $A$ be a unital simple $C^*$-algebra with tracial rank zero
and let $\gamma:\Z_l\to\Aut(A)$ be an action of $\Z_l$.
Suppose that there exists a sequence of projections
$\{e_n\}_{n\in\N}$ satisfying the following property.
\begin{enumerate}
\item For each $i\in\Z_l\setminus\{0\}$,
$\lVert e_n\gamma^i(e_n)\rVert_2\to0$.
\item $\lVert1-\sum_{i\in\Z_l}\gamma^i(e_n)\rVert_2\to0$.
\item For every $a\in A$, we have
$\lVert ae_n-e_na\rVert_2\to0$.
\end{enumerate}
Then the action $\gamma$ has the tracial Rohlin property
in the sense of \textup{\cite{Ph4}}.
\end{prop}
\begin{proof}
The proof goes in a similar fashion to \cite[Section 4]{K}.
By Lemma \ref{replace}, we may assume that
\[ \lim_{n\to\infty}\lVert ae_n-e_na\rVert=0 \]
for every $a\in A$.
Let ${\cal F}$ be a finite subset of $A$ and
let $\varepsilon>0$.
There exist a projection $p\in A$ and
a unital finite dimensional $C^*$-subalgebra $E\subset pAp$
such that the following are satisfied.
\begin{itemize}
\item For every $a\in {\cal F}$,
$\rVert ap-pa\lVert<\varepsilon$.
\item For every $a\in {\cal F}$, there exists $b\in E$
such that $\rVert pap-b\lVert<\varepsilon$.
\item $\lVert1-p\rVert_2<\varepsilon$.
\end{itemize}
Since $\{\gamma^i(e_n)\}_n$ is a central sequence
for every $i\in\Z_l$,
by using the integration argument as in the proof of Lemma \ref{replace},
we may assume that there exists
a projection $x_n^{(i)}\in A\cap E'$
such that $\lVert x_n^{(i)}-\gamma^i(e_n)\rVert\to0$
for every $i\in\Z_l$.
Put\[ y_n=x_n^{(0)}\left(\sum_{i\neq0}x_n^{(i)}\right)x_n^{(0)}. \]
Then $y_n$ is a positive element lying in
\[ D_n=x_n^{(0)}(A\cap E')x_n^{(0)}. \]
Let $\varepsilon_n=\sup_{\tau\in T(A)}\tau(y_n)$.
By the assumption (1), we have $\varepsilon_n\to0$.
Define continuous functions $g_1$, $g_2$ and $g_3$ on $[0,1]$ by
\[ g_j(t)=\begin{cases}
1 & 0\leq t\leq j/4 \\
1-4(t-j/4) & j/4\leq t\leq(j+1)/4 \\
0 & \text{otherwise}. \end{cases} \]
Put $a_{n,j}=g_j(y_n/\sqrt{\varepsilon})$ in $D_n$.
Then it is easy to see that $a_{n,1}a_{n,2}=a_{n,1}$ and
$a_{n,2}a_{n,3}=a_{n,2}$.
Since $D_n$ has real rank zero,
the hereditary subalgebra $\overline{a_{n,2}D_na_{n,2}}$ has
an approximate identity consisting of projections.
Hence we can find a projection
$f_n\in D_n$ such that $a_{n,3}f_n=f_n$ and
$\lVert a_{n,1}f_n-a_{n,1}\rVert<\varepsilon_n$.
Combining
\[ \tau(y_n^{1/2}(x_n^{(0)}-a_{n,1})y_n^{1/2})
\leq\tau(y_n)\leq\varepsilon_n \]
with
\[ \frac{\sqrt{\varepsilon_n}}{4}(x_n^{(0)}-a_{n,1})
\leq y_n^{1/2}(x_n^{(0)}-a_{n,1})y_n^{1/2}, \]
we get $\tau(x_n^{(0)}-a_{n,1})\leq4\sqrt{\varepsilon_n}$
for all $\tau\in T(A)$.
It follows from $\lVert a_{n,1}f_n-a_{n,1}\rVert<\varepsilon_n$
that $\tau(x_n^{(0)}-f_n)\to0$ uniformly for all $\tau\in T(A)$.
Moreover we have
\[ \lVert f_ny_nf_n\rVert=\lVert f_na_{n,3}y_na_{n,3}f_n\rVert
\leq\lVert a_{n,3}y_na_{n,3}\rVert\leq\sqrt{\varepsilon_n}, \]
and so $\lVert f_nx_n^{(i)}f_n\rVert\leq\sqrt{\varepsilon_n}$
for all $i\in\Z_l\setminus\{0\}$.
Therefore, for every $i\in\Z_l\setminus\{0\}$,
\[ f_n\gamma^i(f_n)=f_n\gamma^i(x_n^{(0)})\gamma^i(f_n)
\approx f_nx_n^{(i)}\gamma^i(f_n) \]
converges to zero as $n\to\infty$.
Since $f_n$ commutes with $p\in E$,
$f_np$ is a projection lying in $pAp\cap E'$.
By replacing $f_n$ with $f_np$, we obtain
\[ \lim_{n\to\infty}\sup_{\tau\in T(A)}
\tau(x_n^{(0)}-f_n)<\varepsilon \]
We still have
\[ \lim_{n\to\infty}f_n\gamma^i(f_n)=0 \]
for all $i\in\Z_l\setminus\{0\}$.

As a consequence, by choosing a sufficiently large $n$,
we can find a projection $f_n$ which satisfies the following.
\begin{itemize}
\item For every $a\in {\cal F}$,
$\lVert af_n-f_na\rVert<4\ep.$
\item For every $i\in\Z_l\setminus\{0\}$,
$\lVert f_n\gamma^i(f_n)\rVert<\ep.$
\item For every $\tau\in T(A)$,
\[ \tau\left(1-\sum_{i\in\Z_l}\gamma^i(f_n)\right)
<l\varepsilon. \]
\end{itemize}
Hence $\gamma$ has the tracial Rohlin property.
\end{proof}
\bigskip

Let $\xa$ be a Cantor minimal system and
let $c:X\to\Z_l$ be a continuous map.
Define a homeomorphism $\alpha\times c\in\Homeo(X\times\Z_l)$ by
\[ (\alpha\times c)(x,k)=(\alpha(x),k+c(x)) \]
for all $(x,k)\in X\times\Z_l$.
Namely $\alpha\times c$ is the skew product extension of $\xa$.
Suppose that $\alpha\times c$ is minimal.
Then $C^*(X\times\Z_l,\alpha\times c)$ is a unital simple $A\T$ algebra
with real rank zero.
Define $\gamma\in\Homeo(X\times\Z_l)$ by
$\gamma(x,k)=(x,k+1)$.
Since $\gamma$ commutes with $\alpha\times c$,
it induces an action $\theta$ of $\Z_l$
on $C^*(X\times\Z_l,\alpha\times c)$.
We would like to see that $\theta$ satisfies the hypothesis
of Proposition \ref{trRohlin}.

Let
\[ {\cal P}=\{X(v,k):v\in V,k=1,2,\dots,h(v)\} \]
be a Kakutani-Rohlin partition of $\xa$.
We may assume that the function $c$ is constant
on each clopen set belonging to ${\cal P}$.
For a given $\varepsilon>0$, it is possible to choose ${\cal P}$
so that $h(v)$ is greater than $\varepsilon^{-1}$ for all $v\in V$.
Thus, $\mu(R({\cal P}))$ is less than $\varepsilon$
for all $\mu\in M_\alpha$, where $R({\cal P})$ is the roof set.
For every $v\in V$ and $k=1,2,\dots,h(v)$,
define $c(v,k)\in Z_l$ by $c(v,1)=0$ and
\[ c(v,k)=\sum_{i=1}^{k-1}c(\alpha^{i-1}(x)), \]
where $x$ is a point in $X(v,1)$.
We define a clopen subset $U$ of $X\times\Z_l$ by
\[ U=\bigcup_{v\in V}\bigcup_{k=1}^{h(v)}
X(v,k)\times \{c(v,k)\}, \]
and put $e=1_U\in C^*(X\times\Z_l,\alpha\times c)$.
It is easy to check that $e,\theta(e),\dots,\theta^{l-1}(e)$ are
mutually orthogonal and $\sum_{i\in\Z_l}\theta^i(e)=1$.
Clearly $e$ commutes with elements of $C(X\times\Z_l)$.
Furthermore we have
\[ (u^*eu-e)^2\leq 1_{R({\cal P})\times\Z_l}, \]
where $u$ is the implementing unitary of
$C^*(X\times\Z_l,\alpha\times c)$.
It follows that
\[ \lVert u^*eu-e\rVert_2<\varepsilon. \]
Consequently $\theta$ satisfies
the hypothesis of Proposition \ref{trRohlin},
and so it has the tracial Rohlin property.
Indeed, it can be easily checked that the crossed product
\[ C^*(X\times\Z_l,\alpha\times c)\rtimes_{\theta}\Z_l \]
is stably isomorphic to $C^*\xa$ (see \cite{M1}).
\bigskip

Now we turn to minimal dynamical systems on $X\times\T$.
Let $\xa$ be a Cantor minimal system and
$\phi:X\to\Homeo(\T)$ be a continuous map.
Suppose that $\alpha\times\phi$ is minimal and
non-orientation preserving.
Then $\alpha\times o(\phi)$ is a minimal homeomorphism
on $X\times\Z_2$.
Let $\pi$ be the projection from $X\times\Z_2$ to
the first coordinate.
By \cite[Lemma 8.1, 8.3]{LM2},
$\alpha\times o(\phi)\times\phi\pi$ is a minimal
orientation preserving homeomorphism.
Put
\[ A=C^*(X\times\Z_2\times\T,\alpha\times o(\phi)\times\phi\pi). \]
Then, as in the Cantor case, the shift map
$(x,k,t)\mapsto(x,k+1,t)$ commutes with
the minimal homeomorphism $\alpha\times o(\phi)\times\phi\pi$.
Let us denote the corresponding action of $\Z_2$ on $A$ by $\theta$.
It is not hard to see that $\theta$ globally preserves
the subalgebra $C^*(X\times\Z_2,\alpha\times o(\phi))$.
Therefore, as discussed before,
there exists a projection $e\in C(X\times\Z_2)$ satisfying
the hypothesis of Proposition \ref{trRohlin}.
But, $e$ clearly commutes with elements of $C(X\times\Z_2\times\T)$,
and so we can conclude that $\theta$ on $A$ also satisfies
the hypothesis of Proposition \ref{trRohlin}.

As a direct consequence, we have the following.

\begin{thm}\label{untwist}
Let $\alpha\times\phi$ be a minimal non-orientation preserving
homeomorphism on $X\times\T$. Suppose that
\[ A=C^*(X\times\Z_2\times\T,\alpha\times o(\phi)\times\phi\pi) \]
has tracial rank zero.
Then, the action $\theta$ has the tracial Rohlin property.
In particular, $C^*(X\times\T,\alpha\times\phi)$ also
has tracial rank zero.
\end{thm}
\begin{proof}
This is immediate from Proposition \ref{trRohlin} and
\cite[Theorem 2.7]{Ph4}.
\end{proof}

\begin{cor}\label{untwistc}
Let $\xa$ be a Cantor minimal system and
let $\phi:X\rightarrow\Isom(\T)$ be a continuous map.
If $\alpha\times\phi$ is rigid, then
$C^*(X\times\T,\alpha\times\phi)$ has tracial rank zero.
\end{cor}
\begin{proof}
If $\alpha\times\phi$ is orientation preserving, the result
follows from
\cite[Theorem 5.6]{LM2}.
Suppose that $\alpha\times\phi$ is not orientation preserving.
By \cite[Lemma 8.4]{LM2}, $\alpha\times o(\phi)\times\phi\pi$
is rigid.
It follows from \cite[Theorem 5.6]{LM2} that
\[ A=C^*(X\times\Z_2\times\T,\alpha\times o(\phi)\times\phi\pi) \]
has tracial rank zero.
By Theorem \ref{untwist}, we get the conclusion.
\end{proof}

\section{$C^*$-strongly approximate  flip conjugacy}

Let $\xa$ and $\yb$ be two topological transitive systems.
Let $A=C^*\xa$ and $B=C^*\yb$ be crossed products.
In \cite{Tm}, J. Tomiyama showed that $\xa$ and $\yb$ are flip conjugate
if and only if there exists an isomorphism $\Phi:A\to B$ such that
$\Phi$ maps $j_{\alpha}(C(X))$ onto $j_{\beta}(C(X)).$

The following is an approximate version of Tomiyama's \CA\, flip conjugacy.

\begin{df}\label{6DCKconj}
Let $\xa$ and $\yb$ be two topological transitive systems.
We say that $\xa$ and $\yb$ are
\textit{$C^*$-strongly approximately flip conjugate}
if there exist sequences of isomorphisms $\phi_n:A\to B$, $\psi_n:B\to A$,
$\chi_n:C(X)\to C(Y)$ and $\omega_n:C(Y)\to C(X)$ such that
$[\phi_n]=[\phi_1]$ in $KL(A,B)$, $[\psi_n]=[\psi_1]$ in $KL(B,A)$
$(\phi_n)_{\natural}=(\phi_1)_{\natural}$,
$(\psi_n)_{\natural}=(\psi_1)_{\natural}$,
$\phi_n^{\sharp}=\phi_1^{\sharp}$ and $\psi_n^{\sharp}=\psi_1^{\sharp}$
for all $n\in\N,$
and
\[ \lim_{n\to\infty}\|\phi_n\circ j_{\alpha}(f)
-j_{\beta}\circ\chi_n(f)\|=0 \]
and
\[ \lim_{n\to\infty}\|\psi_n\circ j_{\beta}(g)
-j_{\alpha}\circ\omega_n(g)\|=0 \]
for all $g\in C(X)$ and $g\in C(Y)$.
\end{df}

Let $A$ and $B$ be two unital separable simple \CA\, with real rank zero
and stable rank one and
suppose that there exists an order \hm\,
\[ \tilde\kappa:(K_0(A),K_0(A)_+,[1_A])\to(K_0(B),K_0(B)_+,[1_B]). \]
Let $\rho_A:K_0(A)\to \Aff(T(A))$ be the \hm\, induced
by $\rho_A([p])(\tau)=\tau(p).$
It follows from \cite{BH} that $\rho_A(K_0(A))$ is dense in $\Aff(T(A)).$
Thus $\tilde\kappa$ gives an affine continuous map
$\kappa_{\natural}:\Aff(T(A))\to \Aff(T(B)).$

In the case that $A$ and $B$ are simple and have real rank zero
and stable rank one, in Definition \ref{6DCKconj} above,
if $[\phi_n]=[\phi_1]$ in $KL(A,B),$
then one must have $(\phi_n)_{\natural}=(\phi_1)_{\natural}.$
Moreover, in this case, $K_1(B)=U(B)/U_0(B)=U(B)/CU(B).$
Therefore $\phi_n^{\sharp}=\phi_1^{\sharp}.$
In other words, in Definition \ref{6DCKconj} above,
if both $A$ and $B$ have tracial rank zero,
then one can omit $(\phi_n)_{\natural}=(\phi_1)_{\natural}$ as well as
$\phi_n^{\sharp}=\phi_1^{\sharp}.$
\bigskip

We identify $\T$ with $\R/\Z$ in this section.
Let $\xa$ and $\yb$ be Cantor minimal systems and
let $\phi:X\to\Isom(\T)$ and $\psi:Y\to\Isom(\T)$ be continuous maps.
For the rest of this section we assume that
both $\alpha\times\phi$ and $\beta\times\psi$ are minimal,
but that neither $\alpha\times\phi$ nor $\beta\times\psi$ are
orientation preserving except in Theorem \ref{CCtoKCT}.
We denote the crossed product algebras arising from
$(X\times\T,\alpha\times\phi)$ and $(Y\times\T,\beta\times\psi)$
by $A$ and $B$, respectively.

We identify $K_i(C(X\times\T))$ with $C(X,\Z)$ for $i=0,1$.
By Lemma 2.5 of \cite{LM2}, we know that
$K_0(A)$ is unital order isomorphic to $K^0\xa=\Coker(\id-\alpha^*)$
and that
$K_1(A)$ is isomorphic to
the direct sum of $\Z$ and $\Coker(\id-\alpha_\phi^*)$.
Note that the torsion subgroup of $\Coker(\id-\alpha_\phi^*)$
is isomorphic to $\Z_2$.

In the argument below, we will regard functions of $C(X,\Z)$
as elements of $K_0(A)$ and $K_1(A)$.
When we need to avoid confusion, we denote the equivalence class
of $f\in C(X,\Z)$ in these groups
by $[f]_0$ and $[f]_1$, respectively.

Let $x_0\in X$.
By Proposition 3.3 of \cite{LM2},
we know that $K_0(A_{x_0})$ is unital order isomorphic to $K_0(A)$
and that
$K_1(A_{x_0})$ is isomorphic to
\[ C(X,\Z)/\{f-\alpha_\phi^*(f):f\in C(X,\Z)\text{ and }f(x_0)=0\}. \]
Furthermore, there exists a natural quotient map
from $K_1(A_{x_0})$ to $K_1(A)$ and its kernel is isomorphic to $\Z$.

Define a function $h_\phi\in C(X,\Z)$ by
\[ h_\phi(x)=\begin{cases}1&o(\phi)(\alpha^{-1}(x))=1\\
0&\text{otherwise}.\end{cases} \]
Then $h_\phi$ is a representative of the torsion element in $K_1(A)$.
Thus $2h_\phi$ is zero in $K_1(A)$,
and so $2h_\phi$ belongs to the kernel of the natural quotient map
from $K_1(A_{x_0})$ to $K_1(A)$.
By an easy observation,
we can see that $2h_\phi$ is the generator of the kernel.
Note that
\[ 1_X-\alpha_\phi^*(1_X)=2h_\phi. \]

Let
\[ {\cal P}=\{X(v,k):v\in V,k=1,2,\dots,h(v)\} \]
be a Kakutani-Rohlin partition for $\xa$.
We denote the roof set of ${\cal P}$ by
\[ R({\cal P})=\bigcup_{v\in V}X(v,h(v)). \]
We also write
\[ \widetilde{\cal P}=\{X(v,k):v\in V,k=1,2,\dots,h(v)-1\}
\cup \{R({\cal P})\}. \]
Suppose that ${\cal P}$ is so finer that $o(\phi)$ is constant
on each clopen set belonging to ${\cal P}$.
We define $o(\phi)_v\in\Z_2$ in the same way as in Section 4 of
\cite{M3}. Namely,
\[ o(\phi)_v=o(\phi)(x)+o(\phi)(\alpha(x))+\dots+
o(\phi)(\alpha^{h(v)-1}(x)), \]
where $x$ is a point in $X(v,1)$.

\begin{lem}\label{torsion}
Let $x_0\in X$ and
let
\[ {\cal P}=\{X(v,k):v\in V,k=1,2,\dots,h(v)\} \]
be a Kakutani-Rohlin partition for $\xa$.
Suppose that $x_0$ belongs to $R({\cal P})$ and that
$o(\phi)$ is constant on each clopen set of $\widetilde{\cal P}$.
Then $h_\phi$ is equivalent to
\[ \sum_{o(\phi)_v=1}1_{X(v,h(v))} \]
in $K_1(A_{x_0})$.
\end{lem}

\begin{proof}
We assume that $o(\phi)$ is zero on the roof set.
The other case can be shown similarly.
For every $v\in V$, set
\[ E_v=\{1\leq k\leq h(v):o(\phi)\text{ is not zero on }X(v,k)\}. \]
Then we have
\[ h_\phi=\sum_{n\in V}\sum_{k\in E_v}1_{X(v,k+1)}. \]
Let $k_1<k_2<\dots<k_n$ be the arranged list of elements in $E_v$.
For every $i=1,2,\dots,n$, let $l_i=h(v)-k_i+1$.
It is easy to see that
\[ (\alpha_\phi^*)^{l_i}(1_{X(v,k_i+1)})=(-1)^{n-i}1_{X(v,h(v))}. \]
Hence, if $n$ is even, then
\[ \sum_{i=1}^n(\alpha_\phi^*)^{l_i}(1_{X(v,k_i+1)})=0. \]
If $n$ is odd, then
\[ \sum_{i=1}^n(\alpha_\phi^*)^{l_i}(1_{X(v,k_i+1)})
=1_{X(v,h(v))}. \]
Therefore $h_\phi$ is equivalent to
\[ \sum_{o(\phi)_v=1}1_{X(v,h(v))} \]
in $K_1(A_{x_0})$.
\end{proof}
\bigskip

Let $\sigma$ be an element of the topological full group of $\alpha$.
Then there exists a continuous function $n:X\to\Z$ such that
$\sigma(x)=\alpha^{n(x)}(x)$ for all $x\in X$.
Put $X_k=n^{-1}(k)$ for $k\in\Z$.
Note that $X_k$ is a clopen subset of $X$.
We define an automorphism $\sigma_\phi^*$ on $C(X,\Z)$ by
\[ \sigma_\phi^*(f)=\sum_{k\in\Z}(\alpha_\phi^*)^k(f1_{X_k}) \]
for $f\in C(X,\Z)$.
In other words,
\[ \sigma_\phi^*(f)(x)=(-1)^{c(x)}f(\sigma^{-1}(x)), \]
where $c:X\to\Z_2$ is defined by
\[ c(x)=\begin{cases}
\sum_{i=1}^{n(\sigma^{-1}(x))}o(\phi)(\alpha^{-i}(x))
& n(\sigma^{-1}(x))>0 \\
0 & n(\sigma^{-1}(x))=0 \\
\sum_{i=1}^{-n(\sigma^{-1}(x))}o(\phi)(\alpha^{i-1}(x))
& n(\sigma^{-1}(x))<0. \end{cases} \]
From the $C^*$-algebraic viewpoint,
this definition can be interpreted as follows.
Define $\tilde{\sigma}\in\Homeo(X\times\T)$ by
\[ \tilde{\sigma}(x,t)=(\alpha\times\phi)^{n(x)}(x,t). \]
Thus $\tilde{\sigma}$ belongs to the topological full group
of $\alpha\times\phi$.
Clearly $\tilde{\sigma}$ induces an automorphism of
$K_1(C(X\times\T))$.
Under the identification of $K_1(C(X\times\T))$ with $C(X,\Z)$,
we can see that this automorphism agree with $\sigma_\phi^*$.

\begin{lem}\label{torsionII}
For any $x_0\in X$, $m\in\Z$ and
any nonempty clopen subset $O$ of $X$,
there exists a clopen set $U\subset O$ such that
$1_U$ is equivalent to $mh_\phi$ in $K_1(A_{x_0})$.
Moreover, there exists $\sigma\in[[\alpha]]$ such that
$\sigma(x)=x$ for all $x\in U^c$ and
$\sigma_\phi^*(1_U)=-1_U$.
\end{lem}
\begin{proof}
At first we deal with the case $m=1$.
Since $\alpha\times o(\phi)$ is a minimal homeomorphism
on $X\times\Z_2$, there exists $N\in\N$ such that
\[ \bigcup_{i=0}^N(\alpha\times o(\phi))^i(O\times\{0\})
\supset X\times\{0\}. \]
It follows that, for any $x\in X$, there exists
$i\in\{0,1,\dots,N\}$ such that
\[ (\alpha\times o(\phi))^{-i}(x,0)\in O\times\{0\}. \]
Choose a Kakutani-Rohlin partition
\[ {\cal P}=\{X(v,k):v\in V,k=1,2,\dots,h(v)\} \]
so that the following are satisfied:
\begin{itemize}
\item The roof set $R({\cal P})$ contains $x_0$.
\item $h(v)$ is greater than $N$ for every $v\in V$.
\item $1_O$ and $o(\phi)$ are constant on each clopen set
belonging to $\widetilde{\cal P}$.
\end{itemize}
By the choice of $N$, for each $v\in V$, we can find
$k_v\in\{1,2,\dots,h(v)\}$ such that
$X(v,k_v)$ is contained in $O$ and
$(\alpha_\phi^*)^{h(v)-k_v}(1_{X(v,k_v)})=1_{X(v,h(v))}$.
Put
\[ U=\bigcup_{o(\phi)_v=1}X(v,k_v). \]
Then, by Lemma \ref{torsion},
$1_U$ is equivalent to $h_\phi$ in $K_1(A_{x_0})$.

Let $\sigma$ be the first return map on $U$.
By defining $\sigma(x)=x$ for $x\in U^c$,
we can regard $\sigma$ as an element of $[[\alpha]]$.
We claim $\sigma_\phi^*(1_U)=-1_U$.
There exists $n\in C(X,\Z)$ such that
$\sigma(x)=\alpha^{n(x)}(x)$ for all $x\in X$.
Define $\tilde{\sigma}\in\Homeo(X\times\Z_2)$ by
\[ \tilde{\sigma}(x,k)=(\alpha\times o(\phi))^{n(x)}(x,k). \]
By the choice of $U$, we can see the following.
\begin{itemize}
\item If $x\in X(v,k_v)\subset U$, then
$(\alpha\times o(\phi))^{h(v)-k_v}(x,0)
=(\alpha^{h(v)-k_v}(x),0)$.
\item If $x\in R({\cal P})$, $\alpha(x)\in X(v,1)$ and
$o(\phi)_v=0$, then
$(\alpha\times o(\phi))^{h(v)}(x,0)
=(\alpha^{h(v)}(x),0)$.
\item If $x\in R({\cal P})$, $\alpha(x)\in X(v,1)$ and
$o(\phi)_v=1$, then
$(\alpha\times o(\phi))^{k_v}(x,0)
=(\alpha^{k_v}(x),1)$.
\end{itemize}
It follows that
$\tilde{\sigma}(x,0)=(\sigma(x),1)$ for all $x\in U$.
Hence we have $\sigma_\phi^*(1_U)=-1_U$.

We can prove the case $m=-1$ in a similar fashion.

Let us consider the general case.
Suppose $m>1$. Choose non-empty clopen subsets
$O_1,O_2,\dots,O_m\subset O$ which are mutually disjoint.
By applying the argument above to $O_i$,
we get a clopen set $U_i\subset O_i$ such that
$1_{U_i}$ is equivalent to $h_\phi$ in $K_1(A_{x_0})$.
Moreover, there exists $\sigma_i\in[[\alpha]]$ such that
$\sigma_i(x)=x$ for all $x\in U_i^c$ and
$\sigma_{i\phi}^*(1_{U_i})=-1_{U_i}$.
Let $U=\bigcup_{i=1}^mU_i$ and
$\sigma=\sigma_1\sigma_2\dots\sigma_m$.
Then, $1_U$ is equivalent to $mh_\phi$ in $K_1(A_{x_0})$,
$\sigma(x)=x$ for all $x\in U^c$ and
$\sigma_\phi^*(1_U)=-1_U$.

When $m$ is less than $-1$, a similar proof is valid.
\end{proof}
\bigskip

We would like to show that $C^*$-strongly approximate flip conjugacy
implies approximate $K$-conjugacy under the assumption that
both systems are rigid.
If the systems are rigid, then $A$ and $B$ has tracial rank zero.
Hence, when two isomorphisms from $B$ to $A$ induce the same element
in $KL(B,A)$, we can conclude that
they are approximately unitarily equivalent.
Thus, we may assume that there exist an isomorphism $\Psi:B\to A$,
a sequence of unitaries $w_n\in A$ and
a sequence of isomorphisms $\chi_n:C(Y\times\T)\to C(X\times\T)$
such that
\[ \lim_{n\to\infty}\lVert v_n\Psi(g)v_n^*-\chi_n(g))\rVert=0 \]
for all $g\in C(Y)$.
The isomorphism $\Psi$ induces
a unital order isomorphism $\kappa_0:K_0(B)\to K_0(A)$ and
an isomorphism $\kappa_1:K_1(B)\to K_1(A)$.

Let
\[ {\cal Q}=\{Y(w,l):w\in W,l=1,2,\dots,h(w)\} \]
be a Kakutani-Rohlin partition for $\yb$
such that $o(\psi)$ is constant
on each clopen set of $\widetilde{\cal Q}$.
The above hypothesis implies that
there exists an isomorphism $\chi:C(Y\times\T)\to C(X\times\T)$
such that the following conditions hold (we omit the index $n$
to simplify the notation):
for every $U\in {\cal Q}$,
\[ [\chi_{0*}(1_U)]_0=\kappa_0([1_U]_0) \]
in $K_0(A)$ and
\[ [\chi_{1*}(1_U)]_1=\kappa_1([1_U]_1) \]
in $K_1(A)$.

Keeping these notations,
we will show the approximate $K$-conjugacy.
The proof goes by perturbing $\chi$
with elements of the topological full group $[[\alpha]]$.

\begin{lem}\label{modify}
Let $x_0\in X$ and let
\[ g=\sum_{o(\psi)_w=1}1_{Y(w,h(w))}. \]
Then there exists a homeomorphism $\sigma\in[[\alpha]]$
such that the following conditions are satisfied.
\begin{enumerate}
\item $[\sigma^*\chi_{0*}(1_U)]_0
=\kappa_0([1_U]_0)$ in $K_0(A)$ and
$[\sigma_\phi^*\chi_{1*}(1_U)]_1
=\kappa_1([1_U]_1)$ in $K_1(A)$
for all $U\in {\cal Q}$.
\item For every $U\in \widetilde{\cal Q}\setminus\{R({\cal Q})\}$,
$\sigma_\phi^*\chi_{1*}(1_U)$ is equivalent to
$\sigma_\phi^*\chi_{1*}(\beta_\psi^*(1_U))$ in $K_1(A_{x_0})$.
\item $\sigma_\phi^*\chi_{1*}(g)$ is equivalent to
$h_\phi$ in $K_1(A_{x_0})$.
\end{enumerate}
\end{lem}

\begin{proof}
We have to remark that (1) is automatically satisfied
if we choose $\sigma$ in $[[\alpha]]$.

At first let us consider (3).
Suppose that there exist a homeomorphism $\gamma:Y\to X$ and
a continuous map $\omega:Y\to\Z_2$ such that
$\chi_{0*}$ and $\chi_{1*}$ are given by
\[ \chi_{0*}(f)(x)=f(\gamma^{-1}(x)) \]
and
\[ \chi_{1*}(f)(x)=(-1)^{\omega\gamma^{-1}(x)}f(\gamma^{-1}(x)) \]
for $f\in C(Y,\Z)$ and $x\in X$.
By Lemma \ref{torsion}, $[g]_1$ is the unique torsion element of $K_1(B)$.
Then $\kappa_1([g]_1)$ must be $[h_\phi]_1$, because $\kappa_1$ is
an isomorphism.
We already have
$[\chi_{1*}(g)]_1=\kappa_1([g]_1)$ in $K_1(A)$.
It follows that $\chi_{1*}(g)$ is equivalent to $(2n+1)h_\phi$
in $K_1(A_{x_0})$ for some $n\in\Z$.
Choose $w_0\in W$ and $y_0\in Y$ such that
$o(\psi)_{w_0}=1$ and $y_0\in Y(w_0,h(w_0))$.
We have two possibilities: $\omega(y_0)=0$ or $\omega(y_0)=1$.
We assume $\omega(y_0)=0$.
The other case can be dealt with in a similar fashion.
We can find a clopen neighborhood $O$ of $y_0$ so that
$O\subset Y(w_0,h(w_0))$ and $\omega(y)=0$ for all $y\in O$.
By Lemma \ref{torsionII},
we can find a clopen set $O'\subset \gamma(O)$
such that $1_{O'}$ is equivalent to $nh_\phi$ in $K_1(A_{x_0})$.
Besides, there exists $\sigma\in[[\alpha]]$ such that
$\sigma_\phi^*(1_{O'})=-1_{O'}$ and
$\sigma(x)=x$ for all $x\notin O'$.
Evidently we have
$\sigma_\phi^*\chi_{1*}(1_U)=\chi_{1*}(1_U)$
for all $U\in{\cal Q}\setminus\{Y(w_0,h(w_0))\}$,
because the support of $\sigma$ is contained
in $\gamma(Y(w_0,h(w_0)))$.
When $U=Y(w_0,h(w_0))$, for $x\in X$,
\[ \sigma_\phi^*\chi_{1*}(1_{Y(w_0,h(w_0))})(x)
=(1-21_{O'}(x))(-1)^{\omega(\gamma^{-1}(x))}
1_{Y(w_0,h(w_0))}(\gamma^{-1}\sigma^{-1}(x)). \]
Hence we have
\[ \sigma_\phi^*\chi_{1*}(1_{Y(w_0,h(w_0))})
=\chi_{1*}(1_{Y(w_0,h(w_0))})-21_{O'}. \]
It follows that
\[ \sigma_\phi^*\chi_{1*}(g)
=\chi_{1*}(g)-21_{O'}, \]
and this is equivalent to $(2n+1)h_\phi-2nh_\phi=h_\phi$
in $K_1(A_{x_0})$.
Thus (3) is achieved.

Next, in order to achieve (2), we would like to further
perturb $\sigma_\phi^*\chi_{1*}$ obtained above.
To simplify the notation,
we write $\sigma_\phi^*\chi_{1*}$ obtained above
by $\chi_{1*}$.
Choose $w_0\in W$ arbitrarily.
Put $U=Y(w_0,h(w_0)-1)$.
Since $o(\psi)$ is constant on $U$,
\[ [\chi_{1*}(1_U)]_1
=[\chi_{1*}(\beta_\psi^*(1_U))]_1. \]
Therefore there exists $m\in\Z$ such that
$\chi_{1*}(1_U)+2mh_\phi$ is equivalent to
$\chi_{1*}(\beta_\psi^*(1_U))$ in $K_1(A_{x_0})$.
In a similar fashion to the argument in the preceding paragraph,
we can find $\sigma\in[[\alpha]]$
whose support is contained in $\gamma(U)$ and
$\sigma_\phi^*\chi_{1*}(1_U)$ is equivalent to
$\chi_{1*}(1_U)+2mh_\phi$ in $K_1(A_{x_0})$.
Hence we can conclude that
$\sigma_\phi^*\chi_{1*}(1_U)$ is equivalent to
$\sigma_\phi^*\chi_{1*}(\beta_\psi^*(1_U))
=\chi_{1*}(\beta_\psi^*(1_U))$ in $K_1(A_{x_0})$.

By repeating this procedure,
we can achieve the condition (2) finally.
\end{proof}

The following technical lemma plays a critical role in the proof
of Theorem \ref{CCtoKCT}.

\begin{lem}\label{key}
There exists a homeomorphism $\sigma\in[[\alpha]]$ such that
the following conditions are satisfied.
\begin{enumerate}
\item[\rm(a)] For every $U\in \widetilde{\cal Q}$,
we have
$\alpha^*\sigma^*\chi_{0*}(1_U)
=\sigma^*\chi_{0*}\beta^*(1_U)$.
\item[\rm(b)] For every $U\in \widetilde{\cal Q}$,
we have $\alpha_\phi^*\sigma_\phi^*\chi_{1*}(1_U)
=\sigma_\phi^*\chi_{1*}\beta_\psi^*(1_U)$.
\end{enumerate}
\end{lem}

\begin{proof}
Let $x_0\in X$.
By using Lemma \ref{modify},
we can perturb $\chi:C(Y\times\T)\to C(X\times\T)$
with an element of $[[\alpha]]$
so that the following are satisfied.
\begin{enumerate}
\item $[\chi_{0*}(1_U)]_0=\kappa_0([1_U]_0)$ in $K_0(A)$ and
$[\chi_{1*}(1_U)]_1=\kappa_1([1_U]_1)$ in $K_1(A)$
for all $U\in {\cal Q}$.
\item For every $U\in \widetilde{\cal Q}\setminus\{R({\cal Q})\}$,
$\chi_{1*}(1_U)$ is equivalent to
$\chi_{1*}(\beta_\psi^*(1_U))$ in $K_1(A_{x_0})$.
\item $\chi_{1*}(g)$ is equivalent to
$h_\phi$ in $K_1(A_{x_0})$.
\end{enumerate}
Suppose that there exist a homeomorphism $\gamma:Y\to X$ and
a continuous map $\omega:Y\to\Z_2$ such that
$\chi_{0*}$ and $\chi_{1*}$ are given by
\[ \chi_{0*}(f)(x)=f(\gamma^{-1}(x)) \]
and
\[ \chi_{1*}(f)(x)=(-1)^{\omega\gamma^{-1}(x)}f(\gamma^{-1}(x)) \]
for $f\in C(Y,\Z)$ and $x\in X$.
Let
\[ {\cal P}=\{X(v,k):v\in V,k=1,2,\dots,h(v)\} \]
be a Kakutani-Rohlin partition for $\xa$ such that
the roof set $R({\cal P})$ contains $x_0$.
By choosing ${\cal P}$ sufficiently finer,
we may assume the following.
\begin{itemize}
\item $\omega\gamma^{-1}$ is constant on each clopen set
belonging to ${\cal P}$.
\item $o(\phi)$ is constant on each clopen set
belonging to $\widetilde{\cal P}$.
\item $\gamma^{-1}(X(v,k))$ is contained
in some $Y(w,l)\in{\cal Q}$.
\end{itemize}
Define $\omega'\in C(X,\Z_2)$ as follows.
For $x\in R({\cal P})$, put $\omega'(x)=\omega\gamma^{-1}(x)$.
If $x\in X(v,k)$ and $k\neq h(v)$, then put
\[ \omega'(x)=\omega\gamma^{-1}(x)+\sum_{i=0}^{h(v)-1-k}
o(\phi)(\alpha^i(x)). \]
We remark that $\omega'$ is also constant
on each clopen set of ${\cal P}$.
It is easy to see that
\[ (\alpha_\phi^*)^{h(v)-k}
(\chi_{1*}(1_{\gamma^{-1}(X(v,k))}))
=(-1)^{\omega'(x)}1_{X(v,h(v))}, \]
where $x$ is a point in $X(v,k)$.
For $v\in V$, $w\in W$, $l=1,2,\dots,h(w)$ and $c\in\Z_2$,
let us define a subset $N(v,w,l,c)$ of $\{1,2,\dots,h(v)\}$ by
\[ N(v,w,l,c)=\{k=1,2,\dots,h(v):
X(v,k)\subset \gamma(Y(w,l)),\omega'|_{X(v,k)}=c\}. \]
Then, $\chi_{0*}(1_{Y(w,l)})$ is equivalent to
\[ \sum_{v\in V}(\#N(v,w,l,0)+\#N(v,w,l,1))1_{X(v,h(v))} \]
in $K_0(A)$, and
$\chi_{1*}(1_{Y(w,l)})$ is equivalent to
\[ \sum_{v\in V}(\#N(v,w,l,0)-\#N(v,w,l,1))1_{X(v,h(v))} \]
in $K_1(A_{x_0})$.
Hence, the conditions (1) and (2) above tell us
that, by choosing ${\cal P}$ sufficiently finer,
for every $Y(w,l)\in\widetilde{\cal Q}\setminus\{R({\cal Q})\}$
and $v\in V$,
we have
\[ \#N(v,w,l,0)=\#N(v,w,l+1,0)\text{ and }
\#N(v,w,l,1)=\#N(v,w,l+1,1) \]
if $o(\psi)|_{Y(w,l)}=0$, and
\[ \#N(v,w,l,0)=\#N(v,w,l+1,1)\text{ and }
\#N(v,w,l,1)=\#N(v,w,l+1,0) \]
if $o(\psi)|_{Y(w,l)}=1$.
It follows that there exists a permutation $\pi_v$
on $\{1,2,\dots,h(v)\}$ such that the following holds:
if
\[ \pi_v(k)\in N(v,w,l,c), \]
$k\neq h(v)$ and $l\neq h(w)$, then
\[ \pi_v(k+1)\in N(v,w,l+1,c+o(\psi)|_{Y(w,l)}). \]
Moreover, by (3),
we may also assume that, for every $v\in V$,
\[ \sum_{o(\psi)_w=1}
\#N(v,w,h(w),0)-\#N(v,w,h(w),1)=
\begin{cases}0&o(\phi)_v=0\\1&o(\phi)_v=1.\end{cases} \]
Therefore, we can make the permutation $\pi_v$ so that
the following hold:
\begin{itemize}
\item If $\pi_v(k)\in N(v,w,h(w),c)$ and $k\neq h(v)$,
then
\[ \pi_v(k+1)\in N(v,w',1,c+o(\psi)|_{Y(w,h(w))}) \]
for some $w'\in W$.
\item $\pi_v(h(v))\in N(v,w,h(w),0)$ for some $w\in W$.
\end{itemize}
Notice that the latter condition implies
\[ \omega'|_{X(v,\pi_v(1))}
=o(\psi)|_{R({\cal Q})}+o(\phi)_v \]
for all $v\in V$.

We define $\sigma\in[[\alpha]]$ by
\[ \sigma(x)=\alpha^{k-\pi_v(k)}(x) \]
for $x\in X(v,\pi_v(k))$.
Then one can verify that
\[ \gamma^{-1}\sigma^{-1}\alpha\sigma\gamma(U)=\beta(U) \]
for all $U\in\widetilde{\cal Q}$, which means (a).

Define $\theta:X\to\Z_2$ by
\[ \theta(x)=\omega\gamma^{-1}(x)+\omega'(x)
+\omega'(\sigma^{-1}(x)) \]
for $x\in X$.
Then it is easily verified that
\[ \sigma_\phi^*\chi_{1*}(f)(x)
=(-1)^{\theta(x)}f(\gamma^{-1}\sigma^{-1}(x)) \]
for every $f\in C(Y,\Z)$ and $x\in X$.

Let us check the condition (b).
For every $U\in\widetilde{\cal Q}$ and $x\in X$, we have
\begin{align*}
(\alpha_\phi^*\sigma_\phi^*\chi_{1*})(1_U)(\alpha(x))
&=(-1)^{o(\phi)(x)}(\sigma_\phi^*\chi_{1*})(1_U)(x) \\
&=(-1)^{o(\phi)(x)+\theta(x)}1_U(\gamma^{-1}\sigma^{-1}(x)).
\end{align*}
On the other hand,
\begin{align*}
(\sigma_\phi^*\chi_{1*}\beta_\psi^*)(1_U)(\alpha(x))
&=(-1)^{o(\psi)|_U}(\sigma_\phi^*\chi_{1*})
(1_{\beta(U)})(\alpha(x)) \\
&=(-1)^{o(\psi)|_U+\theta(\alpha(x))}
1_{\beta(U)}(\gamma^{-1}\sigma^{-1}\alpha(x)).
\end{align*}
Since $\gamma^{-1}\sigma^{-1}\alpha\sigma\gamma(U)=\beta(U)$,
$\gamma^{-1}\sigma^{-1}(x)$ belongs to $U$ if and only if
$\gamma^{-1}\sigma^{-1}\alpha(x)$ belongs to $\beta(U)$.
Thus, it suffices to show
\[ o(\phi)(x)+\theta(x)=o(\psi)|_U+\theta(\alpha(x)) \]
for $x\in\sigma\gamma(U)$.

Let $x\in X(v,k)$.
We would like to compute $\theta(\alpha(x))+o(\phi)(x)+\theta(x)$.
Suppose $k\neq h(v)$ and $\pi_v(k)\in N(v,w,l,c)$.
Then
\[ \theta(x)=\omega\gamma^{-1}(x)+\omega'(x)
+\omega'(\alpha^{\pi_v(k)-k}(x))
=\omega\gamma^{-1}(x)+\omega'(x)+c. \]
By the construction of $\pi_v$,
we have
\[ \pi_v(k+1)\in N(v,w,l+1,c+o(\psi)|_{Y(w,l)}) \]
if $l\neq h(w)$, and
\[ \pi_v(k+1)\in N(v,w',1,c+o(\psi)|_{Y(w,h(w))}) \]
for some $w'\in W$ if $l=h(w)$.
In either case, we get
\[ \theta(\alpha(x))=\omega\gamma^{-1}(\alpha(x))
+\omega'(\alpha(x))+c+o(\psi)|_{Y(w,l)}. \]
It follows that
\[ \theta(\alpha(x))+o(\phi)(x)+\theta(x)=o(\psi)|_{Y(w,l)}. \]
If $k=h(v)$, then $\alpha(x)$ belongs to $X(v',1)$
for some $v'\in V$.
Since $\omega\gamma^{-1}(x)=\omega'(x)$, we have
\[ \theta(x)=\omega'|_{X(v,\pi_v(h(v)))}=0. \]
On the other hand,
\[ \theta(\alpha(x))
=o(\phi)_{v'}-o(\phi)|_{R({\cal P})}
+\omega'|_{X(v',\pi_{v'}(1))}. \]
Together with $o(\phi)(x)=o(\phi)|_{R({\cal P})}$, we obtain
\[ \theta(\alpha(x))+o(\phi)(x)+\theta(x)
=o(\phi)_{v'}+\omega'|_{X(v',\pi_{v'}(1))}
=o(\psi)_{R({\cal Q})}. \]
Consequently we have
\[ \theta(\alpha(x))+o(\phi)(x)+\theta(x)=o(\psi)|_U \]
for all $U\in\widetilde{\cal Q}$ and $x\in\sigma\gamma(U)$,
which implies the condition (b).
\end{proof}

\begin{lem}\label{replacing}
Suppose that both $(X\times\T,\alpha\times\phi)$
and $(Y\times\T,\beta\times\psi)$ are rigid.
For any finite subset ${\cal F}$ of $C(Y\times\T)$ and
$\ep>0$, there exist a unitary $v\in A$,
a homeomorphism $\gamma:Y\to X$ and
a continuous map $\omega:Y\to\Isom(\T)$ such that
\[ \lVert v\Psi(f)v^*-f\circ(\gamma\times\omega)^{-1}\rVert<\ep \]
for all $f\in {\cal F}$.
\end{lem}

\begin{proof}
We may assume that ${\cal F}$ is of the form
${\cal F}=\{1_U:U\in{\cal Q}\}\cup\{z\}$, where ${\cal Q}$ is
a clopen partition of $Y$.
There exist a unitary $v_0\in A$ and
a homeomorphism $\gamma\times\rho:Y\times\T\to X\times\T$
such that
\[ \lVert v_0\Psi(f)v_0^*-f\circ(\gamma\times\rho)^{-1}\rVert<1/4 \]
for all $f\in {\cal F}$.
Put $\omega(x)=\lambda^{o(\rho_x)}$, where $\lambda\in\Homeo(\T)$ is
given by $\lambda(t)=-t$ for $t\in\T$.
We can find a unitary $v_1\in A$ such that
\[ v_1v_0\Psi(1_U)v_0^*v_1^*
=1_U\circ(\gamma\times\omega)^{-1}=1_{\gamma(U)} \]
for all $U\in{\cal Q}$.
Note that
\[ [v_1v_0\Psi(z1_U)v_0^*v_1^*]_1
=[z1_U\circ(\gamma\times\omega)^{-1}]_1 \]
in $K_1(A)$ for all $U\in{\cal Q}$.
Since the system is rigid, every invariant measure is the product
of a measure on the Cantor set and the Lebesgue measure on $\T$.
Hence, for every tracial state $\tau\in T(A)$ and
every $n\in\Z\setminus\{0\}$ we have
\[ \tau(v_1v_0\Psi(z^n1_U)v_0^*v_1^*)=0 \]
and
\[ \tau(z^n1_U\circ(\gamma\times\omega)^{-1})=0. \]
The tracial rank of $1_UA1_U$ is zero, and so
$U(A)/CU(A)=U(A)/U_0(A).$
Therefore we can apply Corollary \ref{1C} and
get a unitary $v_U\in 1_UA1_U$ such that
\[ \lVert v_Uv_1v_0\Psi(z1_U)v_0^*v_1^*v_U^*
-z1_U\circ(\gamma\times\omega)^{-1}\rVert<\ep. \]
Let $v_2$ be the direct sum of all the $v_U$'s.
Then $v=v_2v_1v_0$ does the work.
\end{proof}

\begin{lem}\label{arrange}
Suppose that $(X\times\T,\alpha\times\phi)$ and
$(Y\times\T,\beta\times\psi)$ are
$C^*$-strongly approximately conjugate.
Let ${\cal Q}$ be a Kakutani-Rohlin partition of $Y$
such that $o(\psi)$ is constant on each clopen set
of $\widetilde{\cal Q}$.
For any $\ep>0$,
there exist a unitary $v\in A$,
a homeomorphism $\gamma\times\omega:Y\times\T\to X\times\T$
and a continuous function $\xi:Y\to\T$ such that
the following are satisfied.
\begin{enumerate}
\item $\lVert v\Psi(f)v^*
-f\circ(\gamma\times\omega)^{-1}\rVert<\ep$
for all $f\in \{1_U:U\in{\cal Q}\}\cup\{z\}$.
\item $\gamma^{-1}\alpha\gamma(U)=\beta(U)$
for all $U\in\widetilde{\cal Q}$.
\item $\phi_{\gamma(y)}(\omega_y(t))
=\omega_{\gamma^{-1}\alpha\gamma(y)}(\psi_y(t))+\xi(y)$
for all $(y,t)\in Y\times\T$.
\end{enumerate}
\end{lem}

\begin{proof}
We may assume that $\ep$ is less than $1/4$.
Since $(X\times\T,\alpha\times\phi)$ and
$(Y\times\T,\beta\times\psi)$ are
$C^*$-strongly approximately conjugate,
there exist a unitary $v_1$ and
a homeomorphism $\gamma_0:Y\times\T\to X\times\T$ such that
\[ \lVert v_1\Psi(f)v_1^*-f\circ\gamma_0^{-1}\rVert<\ep \]
for all $f\in \{1_U:U\in{\cal Q}\}\cup\{z\}$.
From Lemma \ref{replacing}, we may assume that
$\gamma_0$ arises from a cocycle with values in $\Isom(\T)$.
By Lemma \ref{key}, we obtain an element of
the topological full group $[[\alpha]]$
which satisfies (a) and (b) in Lemma \ref{key}.
There exists a unitary $v_2\in A$
which corresponds to this element.
Let $\sigma_0\in[[\alpha\times\phi]]$ be the homeomorphism
on $X\times\T$ induced by $v_2$.
Then
\[ \lVert v_2v_1\Psi(f)v_1^*v_2^*
-f\circ\gamma_0^{-1}\sigma_0^{-1}\rVert<\ep \]
for all $f\in \{1_U:U\in{\cal Q}\}\cup\{z\}$.
Put $v=v_2v_1$ and $\gamma\times\omega=\sigma_0\gamma_0$.
The condition (a) of Lemma \ref{key} implies (2) directly.
The condition (b) of Lemma \ref{key} implies that,
for every $U\in\widetilde{\cal Q}$,
\[ z1_U\circ(\gamma\times\omega)^{-1}\circ(\alpha\times\phi)^{-1} \]
and
\[ z1_U\circ(\beta\times\psi)^{-1}\circ(\gamma\times\omega)^{-1} \]
have the same $K_1$-class in $K_1(C(X\times\T))$. Since
$\phi,\omega$ and $\psi$ take their values in $\Isom(\T)$, there
must exist a continuous function $\xi:U\to\T$ such that
\[ \phi_{\gamma(y)}(\omega_y(t))
=\omega_{\gamma^{-1}\alpha\gamma(y)}(\psi_y(t))+\xi(y) \]
for all $(y,t)\in U\times\T$.
\end{proof}

\begin{lem}
Suppose that $(X\times\T,\alpha\times\phi)$ is rigid.
For any $\xi\in C(X,\T)$, $\ep>0$ and
a finite subset ${\cal F}\subset C(X\times\T)$,
there exists a unitary $w\in A$ such that
\[ \lVert vj_{\alpha}(f)v^*
-j_{\alpha}(f\circ(\id\times R_\xi)^{-1})\rVert<\ep \]
for all $f\in {\cal F}$.
\end{lem}

\begin{proof}
Consider the monomorphism $\lambda:C(X\times \T)\to A$
defined by $\lambda(f)=j_{\alpha}(f\circ(\id\times R_{\xi})^{-1})$
for $f\in C(X\times \T).$
Then $[j_{\alpha}]=[\lambda]$ in $KL(C(X\times\T),A).$
Since $(X\times\T,\alpha\times\phi)$ is rigid,
every invariant measure has the form
$\mu\times m,$ where $\mu$ is an $\alpha$-invariant measure on $X$ and
$m$ is the Lebesgue measure on $\T$ (see Lemma 4.4 of \cite{LM2}).
Since $m$ is invariant under rotations,
for any $\tau\in T(A),$ we have
\[ \tau(j_{\alpha}(f))=\tau(\lambda(f)) \]
for all $f\in C(X\times \T).$
Again, since $(X\times\T,\alpha\times\phi)$ is rigid,
tracial rank of $A$ is zero.
Therefore $U(A)/CU(A)=K_1(A).$
Thus, by Theorem \ref{ITX},
there is a unitary $v\in A$ such that
\[ \|vj_{\alpha}(f)v^*-\lambda(f)\|<\ep \]
for all $f\in {\cal F}.$ The proof is completed.
\end{proof}

Note that in the following statement, we do not assume that
both systems are non-orientation preserving.

\begin{thm}\label{CCtoKCT}
Suppose that $(X\times\T,\alpha\times\phi)$ and
$(Y\times\T,\beta\times\psi)$ are
$C^*$-strongly approximately flip conjugate and
that both systems are minimal and rigid.
Then there exist an isomorphism $\Psi:B\to A$,
a sequence of unitaries $v_n\in A$ and
a sequence of homeomorphisms $\sigma_n:X\times\T\to Y\times\T$
such that the following conditions are satisfied.
\begin{enumerate}
\item $\lim_{n\to\infty}
\lVert v_n\Psi(f)v_n^*-f\circ\sigma_n\rVert=0$
for all $f\in C(Y\times\T)$.
\item $\lim_{n\to\infty}\lVert f\circ\sigma_n(\alpha\times\phi)\sigma_n^{-1}
-f\circ\beta\times\psi\|=0$ for all $f\in C(Y\times \T)$.
\end{enumerate}
\end{thm}

\begin{proof}
Let $\Psi:B\to A$ be the isomorphism
associated with the $C^*$-strongly approximate conjugacy.
Take $\ep>0$ arbitrarily.
Fix a finite subset ${\cal F}\subset C(Y\times \T).$
Without loss of generality, we may assume
that ${\cal F}=\{1_U:U\in\widetilde{\cal Q}\}\cup\{z\}$
for some Kakutani-Rohlin partition
\[ {\cal Q}=\{Y(w,l):w\in W,1\leq l\leq h(w)\} \]
for $\yb.$
It suffices to show that there exist
a unitary $v\in A$ and
a homeomorphism $\sigma:X\times\T\to Y\times\T$ such that
\[ \lVert v\Psi(f)v^*-f\circ\sigma_n\rVert<\ep \]
and
\[ \lVert f\circ\sigma_n(\alpha\times\phi)^{-1}\sigma_n^{-1}
-f\circ(\beta\times\phi)^{-1}\rVert<\ep \]
for all $f\in\{1_U:U\in\widetilde{\cal Q}\}\cup\{z\}$.

It follows from Lemma 2.4 and 2.5 of \cite{LM2} that,
if the system is orientation preserving, then
$K_1$ of the crossed product is torsion free, and
if the system is not orientation preserving, then
$K_1$ must contain a torsion.
Thus, we have two cases:
both $\alpha\times\phi$ and $\beta\times\psi$ are orientation preserving,
or neither of them are so.

Let us consider the non-orientation preserving case.
By Lemma \ref{arrange}, we can find a unitary $v_1\in A$,
a homeomorphism $\gamma\times\omega:Y\times\T\to X\times\T$
and a continuous function $\xi:Y\to\T$ satisfying the following.
\begin{itemize}
\item $\lVert v_1\Psi(f)v_1^*
-f\circ(\gamma\times\omega)^{-1}\rVert<\ep/2$
for all $f\in \{1_U:U\in{\cal Q}\}\cup\{z\}$.
\item $\gamma^{-1}\alpha\gamma(U)=\beta(U)$
for all $U\in\widetilde{\cal Q}$.
\item $\phi_{\gamma(y)}(\omega_y(t))
=\omega_{\gamma^{-1}\alpha\gamma(y)}(\psi_y(t))+\xi(y)$
for all $(y,t)\in Y\times\T$.
\end{itemize}
By applying Lemma 6.2 of \cite{LM2} to the continuous functions
\[ X\ni x\mapsto (-1)^{o(\phi)(x)}\xi(\gamma^{-1}(x))\in\T \]
and $o(\phi):X\to\Z_2$, we obtain $\eta\in C(X,\T)$ such that
\[ \lvert(-1)^{o(\phi)(x)}\xi(\gamma^{-1}(x))+
\eta(x)-(-1)^{o(\phi)(x)}\eta(\alpha(x))\rvert<\ep. \]
Then we have
\begin{align*}
& \phi_{\gamma(y)}(\omega_y(t)+\eta(\gamma(y)))\\
&=\phi_{\gamma(y)}(\omega_y(t))+(-1)^{o(\phi)(\gamma(y))}\eta(\gamma(y))\\
&=\omega_{\gamma^{-1}\alpha\gamma(y)}(\psi_y(t))+\xi(y)
+(-1)^{o(\phi)(\gamma(y))}\eta(\gamma(y))\\
&\stackrel{\ep}{\approx}\omega_{\gamma^{-1}\alpha\gamma(y)}(\psi_y(t))
+\eta(\alpha(\gamma(y)))
\end{align*}
for all $(y,t)\in Y\times\T$.

Hence, when we put
$\sigma=(\gamma\times\omega)^{-1}(\id\times R_\eta)^{-1}$,
one can check
\[ \lVert f\circ\sigma(\alpha\times\phi)^{-1}\sigma^{-1}
-f\circ(\beta\times\phi)^{-1}\rVert<\ep \]
for all $f\in\{1_U:U\in\widetilde{\cal Q}\}\cup\{z\}$.
By applying the lemma above to $\eta$
and $f\circ(\gamma\times\omega)^{-1}$
for $f\in\{1_U:U\in\widetilde{\cal Q}\}\cup\{z\}$,
we can find a unitary $v_2\in A$ such that
\[ \lVert v_2v_1\Psi(f)v_1^*v_2^*
-f\circ\sigma\rVert<\ep \]
for all $f\in\{1_U:U\in\widetilde{\cal Q}\}\cup\{z\}$.
Thus, we get the unitary $v=v_2v_1$.

We now turn to the orientation preserving case. We may assume that
$\phi$ and $\psi$ take their values in rotations on $\T$. The
isomorphism $\Psi$ induces a unital order isomorphism $\Psi_{*0}$
between $K_0(A)\cong K^0\xa\oplus\Z$ and $K_0(B)\cong
K^0\yb\oplus\Z$. By the definition of $C^*$-strongly approximate
conjugacy, we see that the restriction of $\Psi_{*0}$ on $K^0\yb$
gives a unital order isomorphism from $K^0\yb$ onto $K^0\xa$. We
can identify $K_1(A)$ and $K_1(B)$ with $K^0\xa\oplus\Z$ and
$K^0\yb\oplus\Z$ respectively. Since both $A$ and $B$ have tracial
rank zero, by \cite{L3}, there exists an isomorphism $\Phi:B\to A$
such that $\Phi_{0*}=\Psi_{0*}$ and $\Phi_{1*}=\kappa\oplus\id$.

By \cite[Theorem 5.4]{LM1} or \cite[Theorem 3.4]{M3},
there exists a homeomorphism $\gamma:Y\to X$ such that
$\gamma^{-1}\alpha\gamma(U)=\beta(U)$ for every $U\in{\cal Q}$
and $\kappa([1_U])=[1_{\gamma(U)}]$ for every clopen subset
$U$ of $Y$.
Define a continuous function $\xi:Y\to\T$ by
\[ \xi(y)=\phi_{\gamma(y)}(0)-\psi_y(0) \]
for all $y\in Y$.
By applying \cite[Lemma 6.1]{LM2} to the continuous function
$\xi\circ\gamma^{-1}:X\to\T$,
we obtain $\eta\in C(X,\T)$ such that
\[ \lvert\xi(\gamma^{-1}(x))+\eta(x)-\eta(\alpha(x))\rvert<\ep \]
for all $x\in X$.
Then we have
\[ \eta(\gamma(y))+\phi_{\gamma(y)}(0)
=\eta(\gamma(y))+\xi(y)+\psi_y(0)
\stackrel{\ep}{\approx}\psi_y(0)+\eta(\alpha(\gamma(y))) \]
for all $y\in Y$.
Therefore, when we put $\sigma=(\gamma\times R_\eta)^{-1}$,
one can check
\[ \lVert f\circ\sigma(\alpha\times\phi)^{-1}\sigma^{-1}
-f\circ(\beta\times\phi)^{-1}\rVert<\ep \]
for all $f\in\{1_U:U\in\widetilde{\cal Q}\}\cup\{z\}$.
It is easily verified that
$\Phi_{0*}([1_U])=[1_U\circ\sigma]$ in $K_0(A)$ and
$\Phi_{1*}([z1_U])=[z1_U\circ\sigma]$ in $K_1(A)$
for each clopen subset $U$ of $Y$.
By a similar argument to the proof of Lemma \ref{replacing},
we can find a unitary $v\in A$ such that
\[ \lVert v\Phi(f)v^*-f\circ\sigma\rVert<\ep \]
for finitely many $f\in C(Y\times\T)$,
thereby completing the proof.
\end{proof}

\section{Approximate $K$-conjugacy for minimal rigid systems}

The purpose of this section is to present
a $K$-theoretical condition for which two minimal systems
$(X\times\T,\alpha\times\phi)$ and $(Y\times\T,\beta\times\psi)$
are approximately $K$-conjugate.

We first start with the following definition.

\begin{df}\label{Dconjugate}
Let $\xa$ and $\yb$ be dynamical systems and put
$A=C^*\xa$ and $B=C^*\yb$.
We say that $\xa$ and $\yb$ are approximately $K$-conjugate
if there are homeomorphisms
$\sigma_n:X\to Y$ and $\gamma_n:Y\to X$ such that
\[ \lim_{n\to\infty}\|g\circ\sigma_n\circ\alpha\circ\sigma_n^{-1}
-g\circ\beta\|=0\text{ for all }g\in C(Y), \]
\[ \lim_{n\to\infty}\|f\circ\gamma_n\circ\beta\circ\gamma_n^{-1}
-f\circ\alpha\|=0\text{ for all }f\in C(X), \]
and there are isomorphisms
$\psi_n:B\to A$ and $\phi_n:A\to B$ such that
\[ \lim_{n\to\infty}\|j_{\beta}(f\circ\gamma_n)
-\phi_n(j_{\alpha}(f))\|=0
\text{ and }
\lim_{n\to\infty}\|j_{\alpha}(g\circ\sigma_n)
-\psi_n(j_{\beta}(g))\|=0 \]
for all $f\in C(X)$ and $g\in C(Y).$
Moreover, there exists $\kappa\in KL(A,B)$ and an isomorphism
\[ \tilde\kappa:(K_0(A),K_0(A)_+,[1_A],K_1(A),T(A))
\to(K_0(B),K_0(B)_+,[1_B],K_1(B),T(B)) \]
such that $\kappa$ induces $\tilde\kappa$ on $K_*(A),$
$[\phi_n]=\kappa$ and $[\psi_n]=\kappa^{-1}.$
\end{df}

\begin{rem}
Several remarks about the approximate $K$-conjugacy
are in order.

First, if $\alpha$ and $\beta$ are actually conjugate,
then there exists a homeomorphism $\sigma:X\to Y$ such that
$\sigma\circ\alpha\circ\sigma^{-1}=\beta.$
Define $\Phi(\sum_{-L\le j\le L}f_ju_{\alpha}^j)=
\sum_{-L\le j\le L}f_j\circ\sigma^{-1}u_{\beta}^j$ for $f_j\in C(X),$
$-L\le j\le L.$ It is clear that $\Phi$ gives
an isomorphism from $A$ onto $B.$
Therefore  certainly that $\alpha$ and $\beta$ are conjugate
implies that they are approximately $K$-conjugate.

Second, when $TR(A)=TR(B)=0$ (as the case that we study in this
section), one only needs to require that
$\kappa$ induces an order isomorphism:
$(K_0(A), K_0(A)_+, [1_A], K_1(A))\to(K_1(B), K_0(B)_+, [1_B], K_1(B)).$

Third, if we simply require that
$g\circ\sigma_n\circ\alpha\circ\sigma_n^{-1}\to g\circ\beta$ and
$f\circ\gamma_n\circ\beta\circ\gamma_n^{-1}\to f\circ\alpha$
for all $f\in C(X)$ and $g\in C(Y),$
then $\{\sigma_n\}$ and $\{\gamma_n\}$ may have
no consistent information.
In fact, it was shown in \cite{LM1} that
this requirements are too weaker to be interesting enough in general.
For example, given a projection $p\in C(Y),$
we certainly wish that $[j_{\alpha}(p\circ\sigma_n)]$
eventually gives the same element in $K_0(A).$
These $K$-theoretical consistency on the maps $\sigma_n$ eventually
leads to the above definition.
The reader may notice that when $K_i(A),$ $K_i(B),$
$K_i(C(X)),$ $K_i(C(Y))$ ($i=0,1$)
are torsion free (and $X$ and $Y$ are connected),
Definition \ref{Dconjugate} can be greatly simplified further.
\end{rem}

\begin{thm}\label{KConjT1}
Let $X$ and $Y$ be the Cantor sets and let $(X\times\T,\alpha\times\phi)$
and $(Y\times \T, \beta\times \psi)$ be two minimal rigid systems.
Let $A=C^*(X\times\T,\alpha\times\phi)$
and $B=C^*(Y\times\T,\beta\times\psi).$
Suppose that $\phi_x,\psi_y\in\Isom(\T)$ for each $x\in X$ and $y\in Y.$
Then the following are equivalent.
\begin{enumerate}
\item $(X\times\T,\alpha\times\phi)$ and $(Y\times\T,\beta\times\psi)$
are approximately $K$-conjugate;
\item There exists an isomorphism
$$
\kappa:(K_0(A),K_0(A)_+,[1_A],K_1(A))
\to(K_0(B),K_0(B)_+,[1_B],K_1(B)),
$$
and sequences of isomorphisms $\chi_n:C(X\times\T)\to C(Y\times\T)$
and $\chi_n':C(Y\times\T)\to C(X\times \T)$ such that,
for every finitely generated subgroups $G_i\subset K_i(C(X\times\T))$ and
$F_i\subset K_i(C(Y\times\T)),$
$$
\kappa\circ(j_{\alpha})_*|_{G_i}=(j_{\beta}\circ\chi_n)_*|_{G_i}
\text{ and }
\kappa^{-1}\circ(j_{\beta})_*|_{F_i}=(j_{\alpha}\circ\chi_n')_*|_{F_i}
$$
for $i=0,1$ and all sufficiently large $n;$
\item There exists an isomorphism $\Phi:B\to A$,
sequences of unitaries $\{u_n\}\subset A,$ $\{v_n\}\subset B$ and
sequences of homeomorphisms $\sigma_n:X\times\T\to Y\times\T$ and
$\gamma_n:Y\times\T\to X\times\T$ such that
\[ \lim_{n\to\infty}\lVert u_n\Phi(j_{\beta}(g))u_n^*
-j_{\alpha}(g\circ\sigma_n)\rVert=0 \]
and
\[ \lim_{n\to\infty}\lVert g\circ\sigma_n\circ(\alpha\times\phi)
\circ\sigma_n^{-1}-g\circ(\beta\times\psi)\rVert=0 \]
for all $g\in C(Y\times\T),$ and
\[ \lim_{n\to\infty}\lVert v_n\Phi^{-1}(j_{\alpha}(f))v_n^*
-j_{\beta}(f\circ\gamma_n)\rVert=0 \]
and
\[ \lim_{n\to\infty}\lVert f\circ\gamma_n\circ(\beta\times\psi)
\circ\gamma_n^{-1}-f\circ(\alpha\times \phi)\rVert=0 \]
for all $f\in C(X\times\T).$
\end{enumerate}
\end{thm}

\begin{proof}
(1)$\Rightarrow$(2).
This can be verified directly from Definition \ref{Dconjugate}.

(2)$\Rightarrow$(3).
We first note that, either both $K_1(A)$ and $K_1(B)$ are torsion
free or both has torsion.
By Lemma 2.4 and 2.5 of \cite{LM2},
we note that either both $\alpha\times\phi$ and $\beta\times\psi$ are
orientation preserving or both are non-orientation preserving.

It follows from Corollary \ref{untwistc} that $TR(A)=TR(B)=0.$
It then follows from \cite{L3} that there exists an isomorphism
$\Phi:A\to B$ such that $\Phi$ induces $\kappa.$
Define $\Sigma_n:C(X\times\T)\to B$
by $\Sigma_n(f)=j_{\beta}\circ\chi_n(f).$
Then, by the assumption, one has, for each finitely generated
subgroup $G_i\subset K_i(C(X\times \T))$ ($i=0,1$),
$$
(\Phi\circ j_{\alpha})_*|_{G_i}=(\Sigma_n)_*|_{G_i}
$$
for all sufficiently large $n.$
Let $P_0$ be a set of projections in $C(X\times\T)$
which generates $G_0.$
Thus, for any projection $p\in P_0,$
$[\Phi\circ j_{\alpha}(p)]=[\Sigma_n(p)].$
In particular, for any $\tau\in T(B),$
$$
\tau(\Phi(j_{\alpha}(p))=\tau(\Sigma_n(p)).
$$
Since both systems are rigid and $\phi_x,\psi_y\in\Isom(\T),$
each invariant measure has the form
$\mu\times m,$ where $m$ is the normalized Lebesgue measure on $\T.$
It follows that
$$
\tau(\Phi(j_{\alpha}(f)))=\lim_{n\to\infty}\tau(\Sigma_n(f))
$$
for all $f\in C(X\times\T)$ and $\tau\in T(B).$
Note in the case that $TR(B)=0$ $U(B)/CU(B)=U(B)/U_0(B)=K_1(B).$
Thus, by applying Corollary \ref{ITX2},
we obtain a sequence of unitaries $w_n\in U(B)$ such that
$$
\lim_{n\to\infty}\lVert w_n\Phi(j_{\alpha}(f))w_n^*
-j_{\beta}\circ\chi_n(f)\rVert=0
$$
for all $f\in C(X\times\T).$
Exactly the same argument gives
a sequence of unitaries $v_n\in B$ such that
$$
\lim_{n\to\infty}\lVert v_n\Phi^{-1}(j_{\beta}(g))v_n^*
-j_{\alpha}\circ\chi_n'(f)\rVert=0
$$
for all $g\in C(Y\times\T).$
It follows from Theorem \ref{CCtoKCT} that (3) holds.

(3)$\Rightarrow$(1). This is immediate.

\end{proof}

\begin{rem}\label{RM1}
Consider the case that both $\alpha\times\phi$ and $\beta\times\psi$
are orientation preserving.
It follows from \cite[Lemma 2.4]{LM2} that $K_i(A)=K^0\xa\oplus\Z$ and
$K_i(B)=K^0\yb\oplus\Z$ for $i=0,1$.
Moreover, the embedding $K^0\xa\to K_0(A)$ is an order isomorphism.
Suppose that there exists an isomorphism
\[ \kappa:(K_0(A),K_0(A)_+,[1_A])
\to(K_0(B),K_0(B)_+,[1_B]) \]
such that $\kappa_0\circ(j_{\alpha})_{*0}$ maps
$K_0(C(X\times\T))\cong K_0(C(X))$ onto
$(j_{\beta})_{*0}(K_0(C(Y\times\T))).$
Thus, the restriction of $\kappa_0$ to $K^0\xa\subset K_0(A)$ gives
a unital order isomorphism.
Then, by Theorem 2.6 of \cite{LM1}, one has an isomorphism
$\lambda:C(X)\to C(Y)$ such that $(j_{\beta})_{*0}\circ\lambda_{*0}
=\kappa\circ(j_{\alpha})_{*0}.$
Note that, in the orientation preserving case,
we have $K_1(A)\cong K_0(A)$ and $K_1(B)\cong K_0(B).$
Define $\chi=\lambda\times\id.$
Then it follows that
$\kappa\circ(j_{\alpha})_{*i}=(j_{\beta}\circ \chi)_{*i}$ for $i=0,1.$
One also has
$\kappa^{-1}\circ(j_{\beta})_{*i}=(j_{\alpha}\circ\chi^{-1})_{*i}$
for $i=0,1.$

Thus if $\alpha\times\phi$ and $\beta\times\psi$ are assumed
to preserve the orientation,
Theorem \ref{KConjT1} can be replaced by the following corollary.
\end{rem}

\begin{cor}\label{CCCtoKCJ}
Let $X$ and $Y$ be the Cantor sets and let $(X\times\T,\alpha\times\phi)$
and $(Y\times\T,\beta\times\psi)$ be two minimal rigid systems.
Suppose that $\phi_x,\psi_y$ are rotations for each $x\in X$ and $y\in Y.$
Let $A=C^*(X\times\T,\alpha\times\phi)$ and
$B=C^*(Y\times \T,\beta\times\psi).$
Then the following are equivalent.
\begin{enumerate}
\item $(X\times\T,\alpha\times\phi)$ and
$(Y\times\T,\beta\times\psi)$ are approximately $K$-conjugate;
\item There is an isomorphism
\[ \kappa:(K_0(A),K_0(A)_+,[1_A])\to
(K_0(B),K_0(B)_+,[1_B],) \]
such that
$\kappa_0$ maps $(j_{\alpha})_*(K_0(C(X\times\T)))$
onto $(j_{\beta})_*(K_0(C(Y\times\T)));$
\item There exists an isomorphism $\Phi:B\to A$,
sequences of unitaries $\{u_n\}\subset A,$ $\{v_n\}\subset B$ and
sequences of homeomorphisms $\sigma_n:X\times\T\to Y\times\T$ and
$\gamma_n:Y\times\T\to X\times\T$ such that
\[ \lim_{n\to\infty}\lVert u_n\Phi(j_{\beta}(g))u_n^*
-j_{\alpha}(g\circ\sigma_n)\rVert=0 \]
and
\[ \lim_{n\to\infty}\lVert g\circ\sigma_n\circ(\alpha\times\phi)
\circ\sigma_n^{-1}-g\circ(\beta\times\psi)\rVert=0 \]
for all $g\in C(Y\times\T),$ and
\[ \lim_{n\to\infty}\lVert v_n\Phi^{-1}(j_{\alpha}(f))v_n^*
-j_{\beta}(f\circ\gamma_n)\rVert=0 \]
and
\[ \lim_{n\to\infty}\lVert f\circ\gamma_n\circ(\beta\times\psi)
\circ\gamma_n^{-1}-f\circ(\alpha\times \phi)\rVert=0 \]
for all $f\in C(X\times\T).$
\end{enumerate}
\end{cor}

\begin{rem}\label{FM2}
Let $\xa$ and $\yb$ be two minimal dynamical systems
and let $A=C^*\xa$ and $B=C^*\yb.$
Let $\sigma_n:X\to Y$ and $\gamma_n:Y\to X$ be homeomorphisms
such that
\[ \lim_{n\to\infty}\lVert g\circ\sigma_n\circ\alpha\circ\sigma_n^{-1}
-g\circ\beta\rVert=0 \]
and
\[ \lim_{n\to\infty}\lVert f\circ\gamma_n\circ\beta\circ\gamma_n^{-1}
-f\circ\alpha\rVert=0 \]
for all $f\in C(X)$ and $g\in C(Y)$.
In Definition \ref{Dconjugate},
we required that there exist isomorphisms
which satisfy other requirements.

However, since $A$ and $B$ are nuclear,
as in Proposition 3.2 of \cite{LM1},
there are sequential morphisms $\psi_n:B\to A$ and $\phi_n:A\to B$
such that
\[ \lim_{n\to\infty}\left\lVert\psi_n(\sum_{i=-m}^ng_iu_{\beta}^i)
-\sum_{i=-m}^ng_i\circ\sigma_nu_{\alpha}^i\right\rVert=0 \]
and
\[ \lim_{n\to\infty}\left\lVert\phi_n(\sum_{i=-m}^nf_iu_{\alpha}^i)
-\sum_{i=-m}^nf_i\circ\gamma_nu_{\beta}^i\right\rVert=0, \]
where $f_i\in C(X)$ and $g_i\in C(Y)$.
Unfortunately, in general, $\{\phi_n\}$ and $\{\psi_n\}$ do not give
isomorphisms (not even \hm s).

Suppose that, for any projection $p\in A$ and
any unitary $w\in A,$ we have $[\phi_n(p)]=[\phi_m(p)]$ and
$[\phi_n(w)=[\phi_m(w)]$ for all sufficiently large $n$ and $m.$
Also assume that $\{\phi_n\}$ induces an order isomorphism
$\kappa_0:K_0(A)\to K_0(B)$ and
an isomorphism $\kappa_1:K_1(A)\to K_1(B).$
If we assume that $TR(A)=TR(B)=0,$
then it follows from \cite{L3} that
there is an isomorphism $\Phi:A\to B$ such that $\Phi_{*i}=\kappa_i$
for $i=0,1.$
Suppose also that, for each projection $q\in B$ and
each unitary $v\in B,$ we have $[\psi_n(q)]=[\psi_m(q)]$ and
$[\psi_n(v)]=[\psi_m(v)]$ for all sufficiently large $n$ and $m.$
Then, from $\phi_n(j_{\alpha}(f))-j_{\beta}(f\circ\gamma_n)\to0$ and
$\psi_n(j_{\beta}(g))-j_{\alpha}(g\circ\sigma_n)\to0,$
one sees that, for every finitely generated subgroups
$G_i\subset K_i(C(X\times\T))$ and $F_i\subset K_i(C(Y\times\T)),$
$$
\kappa_i\circ (j_{\alpha})_{*i}|_{G_i}
=(j_{\beta}\circ \gamma_n)_{*i}|_{G_i}\text{ and }
\kappa_i^{-1}\circ (j_{\beta})_{*i}|_{F_i}
=(j_{\alpha}\circ \sigma_n)_{*i}|_{F_i}
$$
for $i=0,1.$
\end{rem}

Therefore we have the following proposition which also explains why
we choose the term approximately $K$-conjugacy.

\begin{prop}\label{7P}
Let $X$ and $Y$ be the Cantor sets and let
$(X\times\T,\alpha\times\phi)$ and $(Y\times\T,\beta\times\psi)$
be two minimal rigid systems. Suppose that $\phi_x$ and $\psi_y$
are in $\Isom(\T)$ for each $x\in X$ and $y\in Y.$ Denote
$A=C^*(X\times\T,\alpha\times\phi)$ and
$B=C^*(Y\times\T,\beta\times\psi).$ Then, $\alpha\times\phi$ and
$\beta\times\psi$ are approximately $K$-conjugate if the following
hold:
\begin{enumerate}
\item There are homeomorphisms $\sigma_n:X\to Y$ and $\gamma_n:Y\to X$
such that
\[ \lim_{n\to\infty}\lVert g\circ\sigma_n\circ\alpha\circ\sigma_n^{-1}
-g\circ\beta\rVert=0 \]
and
\[ \lim_{n\to\infty}\lVert f\circ\gamma_n\circ\beta\circ\gamma_n^{-1}
-f\circ\alpha\rVert=0 \]
for all $f\in C(X)$ and $g\in C(Y)$.
Suppose that $\Phi_n:A\to B$ and $\Psi_n:B\to A$ are
the sequential morphisms induced by $\{\sigma_n\}$ and $\{\gamma_n\}$
as defined in Remark \ref{FM2}.
\item For any projection $p\in A$ and unitary $v\in A,$
$[\Phi_n(p)]=[\Phi_m(p)]$ and $[\Phi_n(v)]=[\Phi_m(v)]$
for all sufficiently large $n$ and $m,$
and $\{\Phi_n\}$ gives a unital order isomorphism
$\kappa_i:K_i(A)\to K_i(B),$ and
\item for any projection $q\in B$ and unitary $w\in B,$
$[\Psi_n(q)]=[\Psi_m(q)]$ and $[\Phi_n(w)]=[\Phi_m(w)]$
for all sufficiently large $n$ and $m,$
and $\{\Psi_n\}$ gives $\kappa_i^{-1}$ $(i=0,1)$.
\end{enumerate}
\end{prop}

\section{Examples}

In this section, we will give two examples.
One example shows that
two minimal systems are approximately $K$-conjugate
but not flip conjugate.
Another example shows that there are minimal systems
whose associated crossed products are isomorphic as
\CA s (and they are weakly approximately conjugate)
but they are not approximately $K$-conjugate.

\begin{exm}
Let $\yb$ be the odometer system of type $5^\infty$.
Let $\xa$ be the Cantor minimal system
described by Figure 2 in \cite[Section 7]{M1}.
Since both $K^0\xa$ and $K^0\yb$ are unital order isomorphic
to $(\Z[1/5],\Z[1/5]_+,1)$, they are strong orbit equivalent.
But, they are not flip conjugate.
Define $c:X\to\Z_2$ by $c(x)=1$ for all $x\in X$.
Then $[c]$ is a nontrivial element of $K^0\xa\otimes\Z_2\cong\Z_2$.
As explained in \cite[Section 7]{M1},
the skew product extension $(X\times\Z_2,\alpha\times c)$ is
a Cantor minimal system and
$K^0(X\times\Z_2,\alpha\times c)$ is also isomorphic to $\Z[1/5]$.
Besides, the canonical inclusion map $K^0\xa$
into $K^0(X\times\Z_2,\alpha\times c)$ is given by $r\mapsto 2r$.
Thus, we have
\[ K^0(X\times\Z_2,\alpha\times c)/K^0\xa\cong\Z_2. \]
Notice that, if we replace $\xa$ with $\yb$,
we obtain exactly the same conclusion.

Let $\xi:X\to\T$ and $\zeta:Y\to\T$ be continuous functions
and put $\phi_x=R_{\xi(x)}\lambda$ and $\psi_y=R_{\zeta(y)}\lambda$
for all $x\in X$ and $y\in Y$,
where $\lambda\in\Homeo(\T)$ is defined by $\lambda(t)=-t$.
Suppose that $\alpha\times\phi$ and $\beta\times\psi$ are
minimal and rigid.
We denote $A=C^*(X\times\T,\alpha\times\phi)$ and
$B=C^*(Y\times\T,\beta\times\psi)$.
It follows from Corollary \ref{untwistc} that
both $A$ and $B$ have tracial rank zero.
By Lemma 2.5 of \cite{LM2},
$K_0(A)$ and $K_0(B)$ are unital order isomorphic to
$K^0\xa\cong K^0\yb$, and
$K_1(A)\cong K_1(B)\cong\Z\oplus\Z_2$.
Therefore $A$ is isomorphic to $B$.
We remark that $[z1_U]$ is nonzero in the $K_1$-group
if and only if $[1_U]$ is not 2-divisible in the $K_0$-group.

Since $K^0\xa$ is unital order isomorphic to $K^0\yb$,
there exists an isomorphism $\rho:C(X)\to C(Y)$
which achieves the order isomorphism $K^0\xa\cong K^0\yb$,
that is, $[1_U]\mapsto[\rho(1_U)]$ gives the order isomorphism.
(See \cite[Theorem 2.6 (3)]{LM1} for example.
Although we only constructed an order isomorphism
from $C(X,\Z)$ to $C(Y,\Z)$ there,
it can be extended to the isomorphism $\rho$.)
Then, we can check the condition (2) of Theorem \ref{KConjT1}
and conclude that $\alpha\times\phi$ and $\beta\times\psi$ are
approximately $K$-conjugate.
But, they are not flip conjugate
because $\alpha$ is not flip conjugate to $\beta$.
\end{exm}

To present examples of two minimal rigid systems whose associated
crossed products are isomorphic but they are not approximately
$K$-conjugate, by applying \ref{KConjT1} and by applying the
classification of unital simple separable amenable \CA s with
tracial rank zero, one only needs to construct two systems whose
$K$-theory of the associated crossed products are unital order
isomorphic but no $\chi_n$ makes the following diagram
\[ \begin{CD}
K_*(A_\alpha) @>\kappa>> K_*(A_\beta) \\
@A(j_\alpha)_*AA @AA(j_\beta)_*A \\
K_*(C(X\times\T)) @>(\chi_n)_*>> K_*(C(Y\times\T))
\end{CD} \]
commute (locally). In the orientation preserving cases, such
examples have been given (\cite[Example 9.2]{LM2}). In what
follows, we construct two non-orientation preserving minimal rigid
systems whose crossed products are isomorphic but they are not
approximately $K$-conjugate. Besides, we construct them so that
they are also weakly approximately conjugate.

\begin{exm}
Let $\theta_1,\theta_2\in(0,1)$ be two irrational numbers
which are linearly independent over $\Q$.
By cutting $\T$ at $n\theta_1$ and $\theta_2+n\theta_1$
for every $n\in\Z$, we get a Cantor set $X$.
Let us denote the $\theta_1$-rotation on $X$ by $\alpha$.
Then $\xa$ is a Cantor minimal system and
$K^0\xa$ is unital order isomorphic to
\[ (\Z+\Z\theta_1+\Z\theta_2,(\Z+\Z\theta_1+\Z\theta_2)_+,1). \]
By cutting $\T$ at $n\theta_2$ and $\theta_1+n\theta_2$
for every $n\in\Z$, we get another Cantor set $Y$.
Let us denote the $\theta_2$-rotation on $Y$ by $\beta$.
Then $\yb$ is also a Cantor minimal system and
$K^0\yb$ is unital order isomorphic to $K^0\xa$.
Hence $\xa$ and $\yb$ are strong orbit equivalent.

Let $U\subset X$ be a clopen subset corresponding to $[0,\theta_1)$.
Define a continuous function $c:X\to\Z_2$
by $c(x)=1$ if and only if $x\in U$.
The skew product extension $(X\times\Z_2,\alpha\times c)$ is
a Cantor minimal system.
By the computation in (2) of \cite[Section 7]{M1},
we see that $K^0(X\times\Z_2,\alpha\times c)$ is isomorphic to
$\Z^5$ and
\[ K^0(X\times\Z_2,\alpha\times c)/K^0\xa\cong\Z_2\oplus\Z^2. \]
Let $\xi:X\to\T$ be a continuous function and
put $\phi_x=R_{\xi(x)}\lambda^{c(x)}$ for all $x\in X$.
Suppose that $\alpha\times\phi$ is minimal and rigid.
Denote $A=C^*(X\times\T,\alpha\times\phi)$.
By Lemma 2.5 of \cite{LM2}, we have $K_0(A)\cong K^0\xa$ and
$K_1(A)\cong \Z\oplus\Z_2\oplus\Z^2$.

Let $V\subset Y$ be a clopen subset corresponding to $[0,\theta_2)$.
By the same way as in the preceding paragraph,
we consider a minimal rigid homeomorphism $\beta\times\psi$ such that
$o(\psi)(y)=1$ if and only if $y\in V$.
We write $B=C^*(Y\times\T,\beta\times\psi)$.
By Lemma 2.5 of \cite{LM2}, we also have $K_0(B)\cong K^0\yb$ and
$K_1(B)\cong \Z\oplus\Z_2\oplus\Z^2$.
It follows from Corollary \ref{untwistc} that
both $A$ and $B$ have tracial rank zero.
Hence $A$ and $B$ are isomorphic.

It can be easily seen that $PS(\alpha)=PS(\beta)
=PS(\alpha\times o(\phi))=PS(\beta\times o(\psi))=\{1\}$,
where $PS(\cdot)$ denotes the set of periodic spectrum.
Therefore, from Corollary 4.10 of \cite{M3},
$\alpha\times\phi$ and $\beta\times\psi$ are
weakly approximately conjugate.

Nevertheless, we would like to show that
$\alpha\times\phi$ and $\beta\times\psi$ are not
approximately $K$-conjugate.
As in Section 6, we identify $K_i(C(X\times\T))$ and $K_i(C(Y\times\T))$
with $C(X,\Z)$ and $C(Y,\Z)$ for each $i=0,1$.
Note that, as explained in Section 6,
$1_U$ is a representative of the unique torsion element of
$\Coker(\id-\alpha_\phi^*)\subset K_1(A)$ and
$1_V$ is a representative of the unique torsion element of
$\Coker(\id-\beta_\psi^*)\subset K_1(B)$.

Assume that $\alpha\times\phi$ and $\beta\times\psi$ are
approximately $K$-conjugate.
We will show a contradiction.
By Theorem \ref{KConjT1},
there exist isomorphisms
\[ \kappa:(K_0(A),K_0(A)_+,[1_A],K_1(A))
\to(K_0(B),K_0(B)_+,[1_B],K_1(B)) \]
and $\chi:C(X\times\T)\to C(Y\times\T)$ such that
\[ \kappa_0\circ(j_\alpha)_{0*}(1_U)=(j_\beta\circ\chi)_{0*}(1_U) \]
and
\[ \kappa_1\circ(j_\alpha)_{1*}(1_U)=(j_\beta\circ\chi)_{1*}(1_U). \]
Since $1_U$ and $1_V$ are representatives of unique torsion elements
in the $K_1$-groups, we must have
\[ \kappa_1\circ(j_\alpha)_{1*}(1_U)=(j_\beta\circ\chi)_{1*}(1_U)
=(j_\beta)_{1*}(1_V), \]
which implies
\[ \chi_{1*}(1_U)-1_V\in\Coker(\id-\beta_\psi^*). \]
It follows that there exists $h:Y\to\Z$ such that
\[ \chi_{1*}(1_U)-1_V=h-\beta_\psi^*(h). \]
Note that $\chi_{1*}(f)-\chi_{0*}(f)$ belongs to $2C(Y,\Z)$
for all $f\in C(X,\Z)$.
Hence
\[ \chi_{0*}(1_U)-1_V\in h-\beta_\psi^*(h)+2C(Y,\Z). \]
It is easy to see that $\beta_\psi^*(g)-\beta^*(g)$ belongs to
$2C(Y,\Z)$ for all $g\in C(Y,\Z)$, and so we get
\[ \chi_{0*}(1_U)-1_V\in h-\beta^*(h)+2C(Y,\Z). \]
On the other hand,
\[ \kappa_0\circ(j_\alpha)_{0*}(1_U)=(j_\beta\circ\chi)_{0*}(1_U) \]
is equal to
\[ \theta_1\in\Z+\Z\theta_1+\Z\theta_2\cong K_0(B), \]
because $\kappa_0$ is a unital order isomorphism.
But, $1_V$ corresponds to $\theta_2$ in $K_0(B)$.
It follows that $\chi_{0*}(1_U)-1_V$ does not belong to $2K_0(B)$.
In other words, there does not exist $h:Y\to\Z$ such that
\[ \chi_{0*}(1_U)-1_V\in h-\beta^*(h)+2C(Y,\Z), \]
which is a contradiction.
\end{exm}

\flushleft{
\textit{Huaxin Lin\\
e-mail: hxlin@noether.uoregon.edu\\
Department of Mathematics\\
University of Oregon\\
Eugene, Oregon 97403\\
U.S.A.\\}
\bigskip
\textit{Hiroki Matui\\
e-mail: matui@math.s.chiba-u.ac.jp \\
Graduate School of Science and Technology,\\
Chiba University,\\
1-33 Yayoi-cho, Inage-ku,\\
Chiba 263-8522,\\
Japan. }}

\end{document}